\documentclass[journal]{IEEEtran}

\usepackage{cite}
\usepackage{todonotes}
\usepackage{graphicx}
\usepackage{caption}

\usepackage{subcaption}
\usepackage{epstopdf}
\usepackage{url}
\usepackage{stfloats}
\usepackage{latexsym}
\usepackage{amsfonts}
\usepackage{amsmath}
\usepackage{amssymb}
\usepackage{array}
\usepackage{textcomp}
\usepackage{pxfonts}
\usepackage{txfonts}
\usepackage{cases}
\usepackage{psfrag}
\usepackage{epsfig}
\usepackage{setspace}
\usepackage{enumitem}
\usepackage{xcolor}
\usepackage{multirow}
\usepackage{multicol}

\newtheorem{Definition}{Definition}
\newtheorem{Remark}{Remark}
\newtheorem{Lemma}{Lemma}
\newtheorem{Proposition}{Proposition}

\newcommand{\R}{\ensuremath{\mathbb{R}}}
\newcommand{\N}{\ensuremath{\mathbb{N}}}
\newcommand{\C}{\ensuremath{\mathbb{C}}}

\begin{document}

%\title{Two-Dimensional Multi-Scale Phase Congruency Parameter Optimization% for Image Features Detection}

%\title{Phase Congruency Optimization to Enhance Image Features Detection}

%\title{Optimal Setting in Image Features Detection Based \\Two Dimensional Multiscale Phase Congruency}

%\title{Optimal and Automatic Tuning of 2D Multiscale Phase Congruency and Image Features Detection}

%\title{Parameter Optimization and Future Open Problems\\ in Image Features Detection via Phase Congruency}

%\title{Parameter Optimization and Future Open Problems\\ in Image Features Detection via Phase Congruency}
%\title{Enhancement of Image Features Detection by Phase Congruency Parameter Optimization with an Example from Multiple Sclerosis MRI }
\title{\Huge Phase Congruency Parameter Optimization for Enhanced Detection of Image Features for both Natural and Medical Applications}

\author{Seyed Mohammad Mahdi Alavi, and
        Yunyan Zhang% <-this % stops a space
        \thanks{S.M.M. Alavi is with the Hotchkiss Brain Institute and Department of Clinical Neurosciences at the Cumming School of Medicine, University of Calgary, Canada. Email: mahdi.alavi@ucalgary.ca.}
      %  \thanks{S. Sharma is with the Hotchkiss Brain Institute and Biomedical Engineering graduate program at the Cumming School of Medicine, University of Calgary, Canada. Email: shrushrita.sharma@ucalgary.ca.}
          % \thanks{G. Pridham is with the Hotchkiss Brain Institute and Department of Clinical Neurosciences at the Cumming School of Medicine, University of Calgary, Canada. Email: glen.pridham@ucalgary.ca.}
\thanks{Y. Zhang is with the Hotchkiss Brain Institute, Department of Clinical Neurosciences, and Department of Radiology at the Cumming School of Medicine, University of Calgary. Email: yunyan.zhang@ucalgary.ca.}% <-this % stops a space
%\thanks{Corresponding  Address: Y. Zhang, Room 183, Heritage Medical Research Building, 3330 Hospital Dr NW, Calgary, AB, Canada, T2N 4N1.}
}

%\markboth{Submitted to IEEE Transactions on Image Processing, January 21, 2017}%

\maketitle

\begin{abstract}
Following the presentation and proof of the hypothesis that image features are particularly perceived at points where the Fourier components are maximally in phase, the concept of phase congruency (PC) is introduced. Subsequently, a two-dimensional multi-scale phase congruency (2D-MSPC) is developed,  which has been an important tool for detecting and evaluation of image features.  However, the 2D-MSPC requires many parameters to be appropriately tuned for optimal image features detection. %This paper in particular addresses the following questions: 
%1) What optimization criteria should be used for parameter tuning of the 2D-MSPC?
%2) How can we formulate the problem mathematically and compute the parameters of the 2D-SMPC optimally and automatically?
%3) How (or by which metric) can the visual convergence be determined? 
In this paper, we defined a criterion for parameter optimization of the 2D-MSPC, which is a function of its maximum and minimum moments. We formulated the problem in various optimal and suboptimal frameworks, and discussed the conditions and features of the suboptimal solutions. The effectiveness of the proposed method was verified through several examples, ranging from natural objects to medical images from patients with a neurological disease, multiple sclerosis. 

%We also described how the proposed method can provide us with a tool to analyze the visual convergence, a condition where an image does not change with adjustment of setting of parameters. 
\end{abstract}


\IEEEpeerreviewmaketitle

\section{Introduction}

%\subsection{Abbreviations}
%\[ \begin{array}{r l}
%
%
%%\mbox{WM:}      & \mbox{White Matter}  \\
%%\mbox{GM:}       & \mbox{Grey Matter}  \\
%%\mbox{CFD}       & \mbox{CerebroSpinal Fluid}\\
%%\mbox{NAWM:}  & \mbox{Normal Appearing White Matter}  \\
%%\mbox{NAGM:}  & \mbox{Normal Appearing Grey Matter}  \\
%\mbox{IFD:}        & \mbox{Image Features Detection}  \\
%\mbox{PC:}        & \mbox{Phase Congruency}  \\
%\mbox{2D-MSPC:} & \mbox{Two-Dimensional Multi-Scale PC}\\
%%\mbox{FLAIR:}   & \mbox{Fluid Attenuated Inversion Recovery}\\
%%\mbox{DIR}       & \mbox{Double Inversion Recovery}\\
%\mbox{MS:}        & \mbox{Multiple Sclerosis}  \\
%\mbox{MRI:}      & \mbox{Magnetic Resonance Imaging} \\
% \end{array}\] 

\subsection{Motivation and Literature Survey}
Image Features Detection (IFD) is an important topic in the image-processing field \cite{Kovesi1996}. It aims at finding image features including lines, edges, Mach bands, corners, and blobs, by using quantitative methods. Many IFD methodologies have been proposed, which can be classified into two broad categories. First, the IFD methods that are based on dimensional metrics. The majority of the published IFD methods lie into this category including dilation-erosion residue operator \cite{Noble1989}, gradient \cite{Canny1983}, weak membrane \cite{Blake1987}, anisotropic diffusion \cite{Perona1990}, and univalue segment assimilating nucleus \cite{Smith1997} -based methods. The main issue of the dimensional IFD methods is that they are too sensitive to image contrasts and spatial magnifications \cite{Kovesi1996}. Second, the IFD methods that are based on non-dimensional metrics. Phase congruency (PC) is such a method \cite{Morrone1986}. References \cite{Kovesi1996}, \cite{Ziou1998}, \cite{Patel2014}, and \cite{Noble1989} review many of the proposed IFD techniques in details.  

The PC-IFD method is important for several reasons \cite{Kovesi1996, Kovesi1999, Kovesi2003}: 1) It is invariant to the image contrast, because it is not based on intensity gradient; 2) No assumption is made in the PC-IFD formulation and computation; 3) It can detect various features; and 4) It uses the image phase information. The latter might be one of the most significant features of the PC-IFD, because phase is shown to be more informative than magnitude in image processing  \cite{Oppenheim1981}. PC is developed based on the hypothesis that image features are optimally perceived at points where the Fourier series components are maximally in phase, meaning that phases of Fourier series€™ components are similar. This hypothesis was first introduced and verified for Mach bands in \cite{Morrone1986}. Subsequently, its effectiveness for detecting other features types has been verified in \cite{Morrone1986, Venkatesh1990, Kovesi1996, Kovesi1999, Kovesi2003}.

It is shown in \cite{Venkatesh1990} and \cite{Morrone1987} that the maxima of the local energy occur at the maxima of PC and vice versa. Thus, in practice, the PC measure is often obtained by computing the normalized local energy, \cite{Morrone1987, Venkatesh1990, Kovesi1999, Kovesi2003}. Moreover, it has been illustrated that the peaks of PC outcome are higher and more distinct when local energy is computed using window-based approaches than the peaks obtained from local energy over the whole signal, \cite{Kovesi1996}. Consequently, Peter Kovesi proposed a multi-scale PC measurement method based on Gabor wavelets and extended it to solve two dimensions (2D), \cite{Kovesi1996, Kovesi1999, Kovesi2003}. Kovesi further modified the PC formulation to overcome noise and ill-conditioning issues, and introduced a weighting function to penalize PC measures at locations where the spread of frequencies is narrow. The proposed 2D multi-scale PC (2D-MSPC) method has been the basis of many IFD studies  in various fields such as history, media, basic sciences, and medicine \cite{Cao2006, Struc2009, Zhang2012, Obara2012, Mouats2015, Rijal2015, Ziaei2014}.

The current computation of 2D-MSPC requires several parameters to be tuned {\em a priori}. In \cite{Kovesi1996, Kovesi1999, Kovesi2003}, several empirical hints have been provided for parameter tuning. But, the selection of parameters either is not explicitly discussed as in \cite{Cao2006, Struc2009, Zhang2012}, or is performed manually based on trail and error \cite{Obara2012, Mouats2015, Rijal2015,  Ziaei2014}. In \cite{Obara2012} and \cite{Mouats2015}, the parameters were manually optimized to increase the visualization of image features. In \cite{Rijal2015} the parameters were changed by trial and error, and the study found that an improvement in IFD occurred with parameters that maximized the signal-to-noise ratio of the image. In \cite{Ziaei2014}, a set of experiments was also done based on the trial and error, to determine the best, fixed values for computing the maximum moment of PC covariance. The authors reported that there was no clear rule on how the 2D-MSPC parameters should be tuned. Collectively, there is a critical need for a standard parameter tuning method for 2D-MSPC. Finding a solution for this problem can enhance our understanding of the concept of PC and thereby promoting its applications.

\subsection{Contribution of This Paper}
We have classified the 2D-MSPC optimization problem as follows:
\begin{itemize}
\item[Q1:] What optimization criteria should be used for the tuning of 2D-MSPC parameters?
\item[Q2:] How can we formulate the problem mathematically and compute the parameters of 2D-SMPC optimally and automatically?
%\item[Q3:] How can the visual convergence be determined? Visual convergence refers to the condition when some parameters of the 2D-MSPC converge to their lower or upper limits, and the image does not change visually by further increasing or decreasing the parameters. 
\end{itemize}

This paper aims to address these questions. The 2D-MSPC, its definition, computational method and a list of its fundamental parameters are briefly described in section \ref{sec:2dmspc}. The importance of parameter setting in 2D-MSPC IFD, which leads to Q1 and Q2, is discussed in section \ref{sec:pardisc}. In section \ref{sec:opt}, we define a criterion for parameter optimization, based on the maximum and minimum moments of the 2D-MSPC. We formulate the problem in various optimal frameworks, and then describe about a suboptimal solution accordingly. %A feature of the suboptimal solution is that it does not require the computation of the 2D-MSPC momentums at each iteration of the optimization. 
%In section \ref{sec:visconv}, we describe how the method can be used for the assessment of the visual convergence. 
Finally, the effectiveness of the proposed method is verified through several examples from both natural and medical images in section \ref{sec:ex}. This includes two examples from magnetic resonance imaging  (MRI) of patients with multiple sclerosis (MS). We illustrate that the proposed 2D-MSPC optimization is useful to MS lesions detection. %Interested readers are referred to \cite{Cabezas2014, Schmidt2012, Lorenzo2013} for a more comprehensive information on the MS MRIs signal processing.  

\subsection{Nomenclatures}
In this paper the following general notations are used. $\N$, $\R$ and $\C$ denote the integer, real and complex domains, respectively. $\|A\|$ denotes the norm of matrix $A$, where no subscript means that any norm can be used. $\|A\|_F$, $\|A\|_p$,  $\|A\|_\infty$, $\|A\|_2$ and $\|A\|_1$ represent Frobenius-norm ($F-$norm), Schatten $p-$norm, $\infty-$norm, 2-norm and 1-norm of matrix $A$. The consistent norm of matrix $A$ is denoted by $\|A\|_c$. The determinant of matrix $A$ is denoted by $\mbox{det}(A)$. The optimal value of the parameter $a$ is denoted by $a^*$ (a superscript asterisk). It should not be confused with the convolution operator $*$, which is used in this paper as $A*B$.  Note that there are several other specific notations also defined in the text, which will be introduced when they appear.

%The paper is organized as follows. In section \ref{sec:2dmspc}, 2D-MSPC, its definition, computation method and a list of its fundamental parameters are briefly described. The importance of the  parameters values, Q1, Q2 and Q3 are discussed in section \ref{sec:pardisc}. The parameter optimization methods are proposed in section \ref{sec:opt}. This section also provides the answers to Q1 and Q2. In section \ref{sec:visconv}, the visualization convergence  Q3 is addressed. Examples are provided in section \ref{sec:ex}. %The paper ends up with two open problems. 

%>>>>>>>>>>>>>>>>>>>>>>>>>>>>>>>>>>>>>>>>>>>
\section{Two-Dimensional Multi-Scale Phase Congruency}
\label{sec:2dmspc}
\begin{figure}[h]
\centering
\includegraphics[scale=.45]{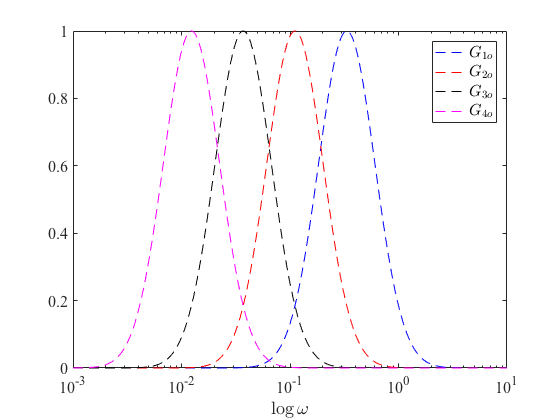}
\caption{Spectra of 4-scale Gabor wavelets.}
\label{fig:gf}
\end{figure}

In this section, fundamental concepts of the 2D-MSPC are briefly described. Interested readers are directed to \cite{Kovesi1996}, \cite{Kovesi1999} and \cite{Kovesi2003} for details. The 2D-MSPC is a combination of one-dimensional PC calculated over several orientations. The 2D-MSPC based on the measurement of local energy is given by 
 \begin{align}
\label{eq:pc}
& \mbox{PC}(x)= \frac{ \sum_{o}  W_o(x) \left \lfloor  E_o(x)-T_o \right \rfloor}{ \sum_{o} \sum_{n}A_{no}(x) + \varepsilon}
\\ \nonumber & o=1,2,\ldots, O
\\ \nonumber & n=1,2,\ldots, N_o
\end{align}
where, $o$ and $n$ denote the orientation and scale indexes, respectively. $O$ is the total number of orientations. $N_o$ is the total number of scales at the $o-$th orientation. $\lfloor \rfloor$ is the floor operator. Note that $\lfloor y \rfloor = y$ if $y>0$, otherwise zero. $E_o$ is the local energy at the $o-$th orientation, which is calculated by using the Hilbert transform as follows
\begin{align}
\label{eq:e}
E_o(x)=\sqrt{F_o^2(x)+H_o^2(x)}
\end{align} 
where, $F_o(x)$ and $H_o(x)$ are the AC component and Hilbert transform of the image signal $I(x)$ at the $o-$th orientation. However, because the Hilbert transform operator is an improper integral, $F_o(x)$ and $H_o(x)$ are computed by convolving the image signal with a pair of even and odd wavelets filters in quadrature as follows \cite{Morrone1987, Venkatesh1990}: 
\begin{align}
\label{eq:fs}
&F_o(x)= \sum_{n} I(x) * M_{no}^{e}\\
\label{eq:hs}
&H_o(s)= \sum_{n} I(x) * M_{no}^{o}%\\
%\nonumber 
%& n=1,2,\ldots,N_o
\end{align}
where, $*$ denotes the convolution operator. $M_{no}^e$ and $M_{no}^o$ are even and odd filters at the $o-$th orientation and $n-$th scale, generated by the logarithmic Gabor function
\begin{align}
\label{eq:lgabor1}
G_{no}(f)=\exp \left(\frac{-\left(\log\frac{f}{\hat{f}_{no}}\right)^2}{2\left( \log \sigma_{no}\right)^2}\right)%~n=1,2,\ldots,N_o
\end{align}
where, $\hat{f}_{no}$ and $\sigma_{no}$ denote the centre frequency and bandwidth of the $n-$th Gabor filer at the $o-$th iteration. The filters are scaled as follows
\begin{align}
\label{eq:lgabor2}
\hat{f}_{no}=\frac{1}{\lambda_{\min o} \times \eta_o^{(n-1)}}%,~n=1,2,\ldots,N_o
\end{align}
where, $\lambda_{\min o}$ is the minimum wavelength (maximum centre frequency) of wavelets in banks of $\{M_{1o}^e, \ldots, M_{No}^e\}$ and $\{M_{1o}^o, \ldots, M_{No}^o\}$,  and $\eta_o$ represents the distance between successive filters in the bank. It is noted that $M_{no}^e$ and $M_{no}^o$ are identical with a $90^{\circ}$ shift. Figure \ref{fig:gf} shows the spectra of 4-scale logarithmic Gabor wavelets at the $o-$th orientation with $\lambda_{\min o}=3$, $\eta_o=3$ and $\sigma_{no}=0.55$, $n=1,\ldots,4$. Note that in the logarithmic frequency scale, the spectra of all Gabor functions are identical. 

Subsequently, the normalization factor in the 2D-MSPC \eqref{eq:pc} is computed in the wavelet's framework as follows
\begin{align}
\label{eq:nf}
\sum_{n}A_{no}=\sum_{n} \sqrt{  (I(x) * M_{no}^{e})^2+(I(x) * M_{no}^{o})^2}.
\end{align}

\begin{table*}[b]
\caption{A list of 2D-MSPC parameters and the existing hints for manual tuning as described in \cite{Kovesi1996, Kovesi1999, Kovesi2003, Kovesi_url}.}
   \centering
     \begin{tabular}{ l  l  l   }
\hline
		%\multicolumn{2}{c}{Parameters} &  Hints  \cite{Kovesi1996} \\ \hline
       \multirow{1}{*}{} & \small Parameter &  \small Tuning hints\\ \hline
       \multirow{2}{*}{\small Weighting Function}&        $c_o$  &   {\scriptsize Between 0 and 1, typical 0.4 or 0.55}      \\ 
         								& $g_o$  & {\scriptsize Typical value, 10}      \\  \hline                  
       \multirow{3}{*}{\small Bank of Filters} & $ \lambda_{\min o}$ & {\scriptsize The smallest value is the Nyquist wavelength of 2 pixels. Because of aliasing 3 pixels or above is suggested.} \\% \cline{2-2}
       								&$\sigma_{no}$         & {\scriptsize The smaller $\sigma_{no}$, the larger the bandwidth of the filter. The following sets are suggested: $[\sigma_{no}=0.85,\eta_o=1.3]$, or }   \\ %\cline{2-2} 
                 						&      $\eta_o$            & {\scriptsize $[\sigma_{no}=0.75,\eta_o=1.6]$ (Bandwidth $\approx$ 1 octave), or $[\sigma_{no}=0.65,\eta_o=2.1]$, or $[\sigma_{no}=0.55,\eta_o=3]$ (Bandwidth $\approx$ 2 octaves)}   \\ \hline
						
	\multirow{1}{*}{\small Noise scaling factor} & $k_o$ &  {\scriptsize Typical value 2, up to 10 or 20 for noisy images} \\ \hline
		 
	\multirow{1}{*}{\small Singularity avoidance} & $\varepsilon$ & {\scriptsize Small scalar, typically 0.0001  } \\ \hline
		 
	\multirow{1}{*}{\small No. of scales} & $N_o$ & {\scriptsize Try values 3 to 6 }\\ \hline
		 
	\multirow{1}{*}{\small No. of orientations} & $O$ &  {\scriptsize A filter spacing of $30^{\circ}$ has found good. } \\ \hline
		 
\end{tabular}
\label{table:par}
\end{table*}

As discussed in \cite{Kovesi1996}, the calculation of PC makes sense only in locations where the spread of frequencies is significant. In 2D-MSPC formulation, $W_o(x)$ is a weighting function, which penalizes the information at locations where the spread of frequencies is narrow at the $o-$th orientation. It is defined by a sigmoid function 
\begin{align}
W_o(x)=\frac{1}{1+\exp \left(g_o(c_o-s_o(x))\right)}
\end{align}
where, $c_o$ denotes the cutoff (mid-point) point of the sigmoid function; $g_o$ is a gain that controls the rate of weighting; $s_o$ is a measure of frequency spread, which ranges between 0 and 1, and is given by
\begin{align}
s_o(x)=\frac{1}{N_o}\frac{\textstyle \sum_{n}A_{no}(x)}{A_{\max}(x)+\varepsilon}.
\end{align}
where,  $\textstyle \sum_{n}A_{no}(x)$ is computed through \eqref{eq:nf} and $A_{\max}(x)$ is the amplitude of the filter pair having maximum response at $x$. 

The parameter $T_o$ is an estimation of noise power that is subtracted from the energy of the signal. It is computed based on the assumptions that the image noise is additive, that the noise power spectrum is constant, and that the image features occur at isolated locations:
\begin{align}
T_o=\mu_{Ro}+k_o \sigma_{Ro}
\end{align}
where, $\mu_{Ro}$ and $\sigma_{Ro}^2$ are the mean and variance of the Rayleigh that describes the noise energy response; $k_o$ is a scaling factor used to estimate the maximum degree of the noise response. %There is little basis to build a more formal model as mentioned in  \cite{Kovesi1996} and references therein.  

When the spread of frequencies is narrow, $E_o$ and $\sum_n A_{no}$ become very small, making the computations to become ill-conditioned. The parameter $\varepsilon$ is to address this issues.  

%\begin{Remark}
%The components and parameters, discussed above, arefundamental. In \cite{Kovesi_url}, other components (for instance a filter to upset the normalization process) have also been added which are neglected from this paper.  Adding new parameters to the proposed optimization methodology is straightforward.  \hfill{$\square$}
%\end{Remark}

Table \ref{table:par} summarizes the parameters discussed above which are fundamental and need to be fine-tuned in the 2D-MSPC calculations, including some hints from  \cite{Kovesi1996, Kovesi1999, Kovesi2003, Kovesi_url}, mainly for manual tuning.

Hereafter, we denote all parameters in a vector format,
\begin{align}
\label{parvec}
\upsilon_o= \left[ \begin{array}{ccccccccc}
    c_o &  g_o & \lambda_{\min o} & \sigma_{no} & \eta_o & k_o & \varepsilon & N_o & O \\ 
  \end{array}\right]
\end{align}
where, $\upsilon_o$ denote the parameter vector at orientation $o$. 

As usual in the numerical optimization, we assume that upper and lower limits of the parameter vector are known {\em a priori}. The upper and lower limits of $\upsilon$ in the vector format are denoted by $\overline{\upsilon}_o$ and $\underline{\upsilon}_o$ respectively, such that
\begin{align}
\underline{\upsilon}_o \leq \upsilon_o \leq \overline{\upsilon}_o
\end{align}
%The upper and lower limits can vary for each orientation in general. 
 
In this paper, $\upsilon_o^*$ denotes the optimal parameter vector and $\mbox{PC}_o^*$ denotes the PC matrix computed with the optimal parameter vector $\upsilon_o^*$, both at orientation $o$.

%\begin{Remark}
%In this study our focus is on the bank of filters' and weighting function's parameter optimization. We set other parameter to fixed values as $k_o=2$, $\varepsilon=0.0001$, $N_o=4$ and $O=6$. \hfill{$\square$}
%\end{Remark}

%Hereafter, we define the parameter vector as:
%\begin{align}
%\theta=\left[
%\begin{array}{c c c c c c c c c}
%  O& N_o  &  & & & & & & \\
%\end{array}
%\right]
%\end{align}

%>>>>>>>>>>>>>>>>>>>>>>>>>>>>>>>>>>>>>>>>>>>
\section{Description of Q1 and Q2}
\label{sec:pardisc}
\begin{figure}[h]
      \centering  \begin{subfigure}[b]{0.25\textwidth}
        		\centering 
                \includegraphics[width=\linewidth]{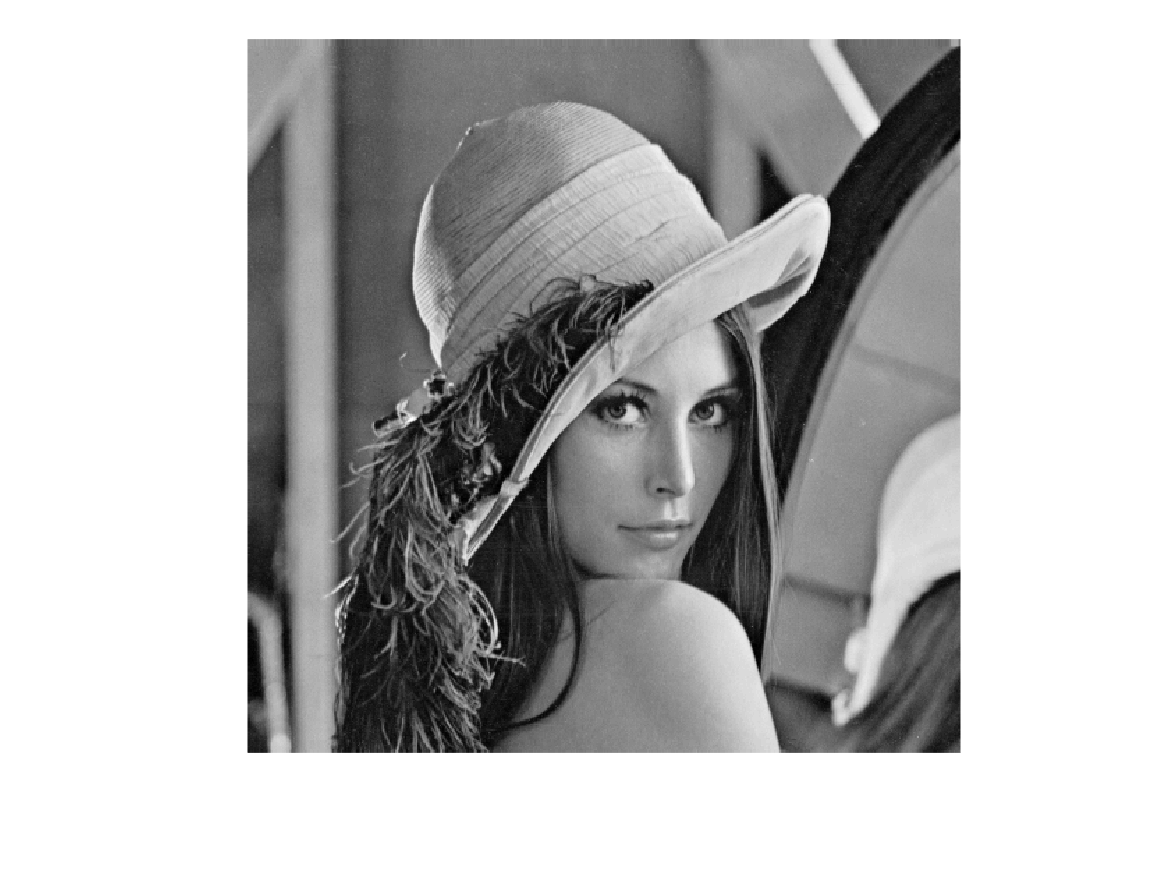}
                \caption{Original image}
                \label{fig:lena-orig}
        \end{subfigure}\\
        \begin{subfigure}[b]{0.25\textwidth}
                \includegraphics[width=\linewidth]{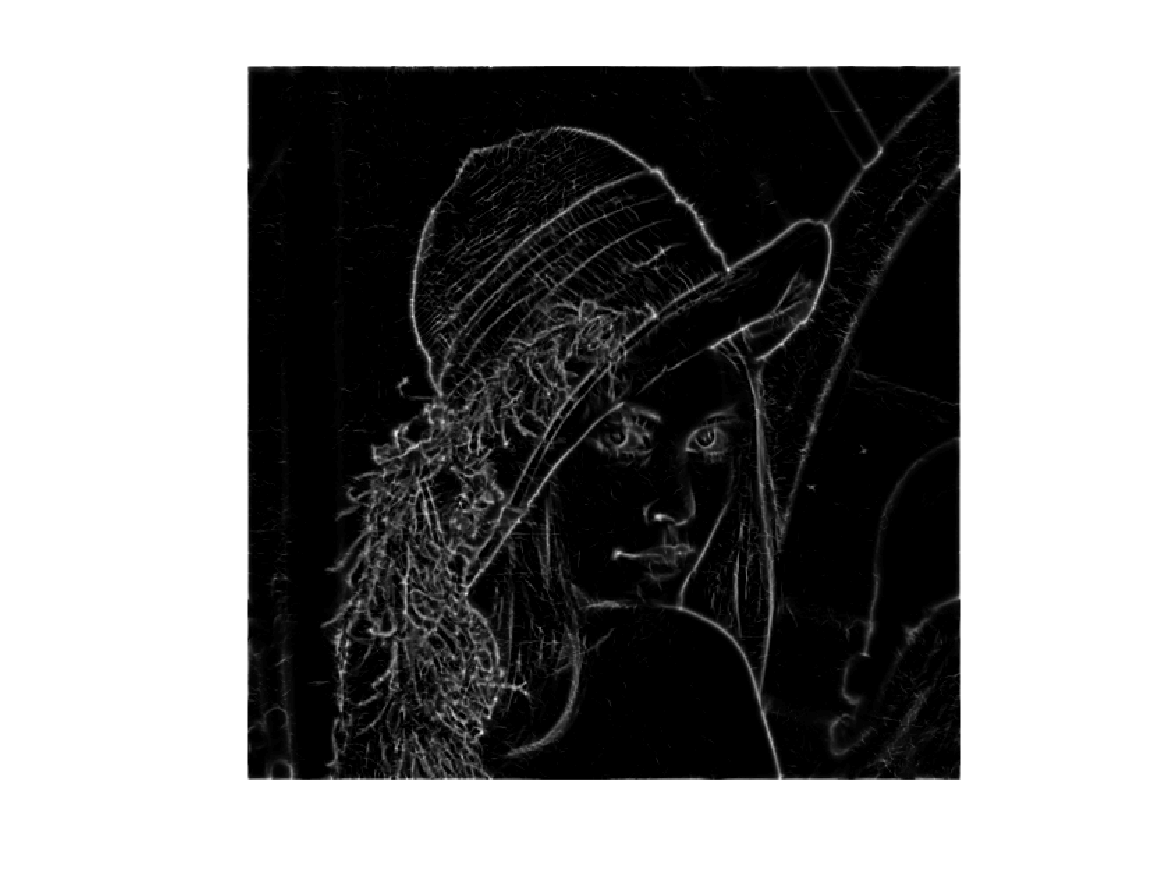}
                \caption{$c_o=0.55$}
                \label{fig:lena-co-p55}
        \end{subfigure}%
        \begin{subfigure}[b]{0.25\textwidth}
                \includegraphics[width=\linewidth]{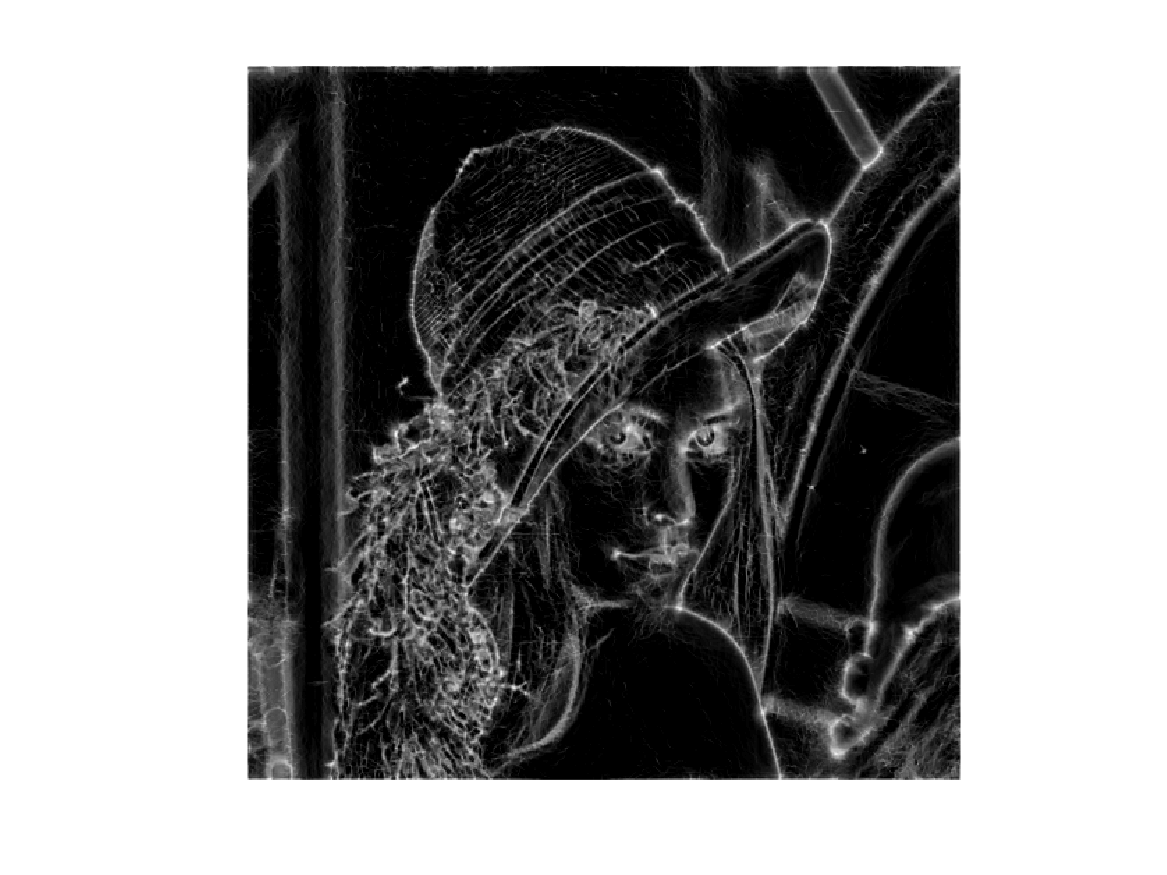}
                \caption{$c_o=0.1$}
                \label{fig:lena-co-p10}
        \end{subfigure}%
%        \begin{subfigure}[b]{0.25\textwidth}
%                \includegraphics[width=\linewidth]{lena-pc-co-p1}
%                \caption{$c_o=0.01$}
%                \label{fig:lena-co-p1}
%        \end{subfigure}
        \caption{Visualization of  a lena image by using the 2D-SMPC, with different values of the cutoff point $c_o$, while other parameters are fixed.}
         \label{fig:lena}
\end{figure}

It is illustrated that visualization of the image will significantly vary by changing only one parameter of the 2D-MSPC \cite{Mouats2015} and \cite{Rijal2015}. Figure \ref{fig:lena} shows the visualization of a lena's image by using 2D-MSPC, with different values of the cutoff point $c_o$. Based on information from Table \ref{table:par}, other parameters of the 2D-MSPC are fixed to: $k_o=2$; $\varepsilon=0.0001$; $N_o=4$; $O=6$; $\lambda_{\min o}=3$; $\eta_o=2.1$; and $g_o=10$, for  $o=1,\dots, O$ and $n=1,\ldots,N_o$. Figure \ref{fig:lena-orig} shows the lena's grayscale image of lena that was used in this simulation. Figures \ref{fig:lena-co-p55} and \ref{fig:lena-co-p10} show the 2D-MSPC images for $c_o=0.55$ and $c_o=0.10$, respectively.  It is seen that by decreasing the cutoff value, image features become increasingly  detectable. There are several features that could not be detected with $c_o=0.55$ (Figure \ref{fig:lena-co-p55}). This issue gives rise to questions Q1 and Q2.

In the next few sections, these questions are addressed.

%>>>>>>>>>>>>>>>>>>>>>>>>>>>>>>>>>>>>>>>>>>>
\section{Parameter Optimization}
\label{sec:opt}
In this paper, several optimal frameworks are proposed for parameter tuning for the 2D-MSPC, which are based on the maximization of the PC momentums. The maximum and minimum moments of the 2D-MSPC, $M$ and $m$ respectively, are given by \cite{Kovesi2003}:
\begin{align}
\label{maxmo}
& M=\frac{1}{2}\left(\alpha+\gamma+\sqrt{\beta^2+(\alpha-\gamma)^2}\right)\\
\label{minmo}
& m=\frac{1}{2}\left(\alpha+\gamma-\sqrt{\beta^2+(\alpha-\gamma)^2}\right)
\end{align}
with
\begin{align}
\label{alpha}
& \alpha=\textstyle \sum_{o}\left(\mbox{PC}_o \cos(\theta_o)\right)^2\\
\label{gamma}
& \gamma=\textstyle \sum_{o}\left(\mbox{PC}_o \sin(\theta_o)\right)^2\\
\label{beta}
& \beta=2 \textstyle \sum_{o}\left(\mbox{PC}_o \cos(\theta_o)\right)\left(\mbox{PC}_o \sin(\theta_o)\right).
\end{align}
where, $\mbox{PC}_o$ and $\theta_o$ represent the PC and axis angle at orientation $o$.

The calculation of momentums corresponds to performing a singular value decomposition to the PC covariance matrix, and thus momentums correspond to the singular values, \cite{Kovesi2003}. On the other hand, singular values indicate the level of increase in energy that can occur between the input and output of a given system. The larger the maximum singular value of a system, the greater the increase in energy, (see Definition 40.2 and descriptions in \S 40.2.1, p. 652,\cite{Levine1995}). Therefore, it may be deduced that momentums represent the maximum of local energies. It was also shown that the maxima of local energy are indications of features in an image, \cite{Venkatesh1990}. The answer to Q1 can now be simplified as such that a potential criterion for parameter tuning of 2D-MSPC should be given in terms of $M$ and $m$.

%The magnitude of the maximum moment determines the feature's significant, the larger the value, the more the feature's significance. The minimum moment gives also insights about the orientation of the feature. The large minimum moment points out that the feature point has also a strong 2D component.
From the potential relationships between momentums, singular values, local energies, 2D-MSPC and image features, it may be deduced that the image features detection can be enhanced by increasing the maximum and minimum moments of the PC. In other words, optimal values of 2D-MSPC parameters are those that maximize $M$ and $m$. Both $M$ and $m$ are matrices and their size is equal to the size of $\mbox{PC}_o$. Thus, the optimal parameters of 2D-MSPC should be obtained by solving a matrix optimization problem (an answer to Q2).

In order to formulate the problem, we define a cost function $\mathcal{M}$:%, which is a linear combination of the 2D-MSPC momentums as follows:
\begin{align}
\label{Mm}
\mathcal{M}=\mu_1 \times  M+\mu_2 \times  m
\end{align}
where, $\mu_i \in [0,1]$, $i=1,2$. %If $\mu_1=1$ and $\mu_2=0$, the maximum moment is going to be maximized. With $\mu_1=0$ and $\mu_2=1$, the minimum moment $m$ is maximized. A combination of $M$ and $m$ are maximized with $\mu_i \neq 0$, $i=1,2$. 
%It is clear that maximization of $M$ and $m$ is equivalent to the maximization of $\mathcal{M}$. 
The cost function \eqref{Mm} enables us to not only maximize both $M$ and $m$, but also to manage their interactions.

As mentioned above, the question of parameter tuning for 2D-MSPC is a matrix optimization problem. Matrix maximization has widely been used in optimal (experimental) designs  \cite{Fredov1972, Goodwin1977, Atkinson1982, Pronzato1985, Berger1994, Pukelsheim1993}, where the Fisher information matrix \cite{Fisher1971} is maximized in order to reduce the Cram\' er-Rao bound \cite{Cramer1946, Rao1945}. The Fisher information matrix is maximized by maximizing some real-valued summary statistics (optimal) criteria), \cite{Atkinson2007}. The majority of matrix optimization methods are based on the maximization of the determinant of a matrix, referred to as D-optimal problem \cite{Berger1994, StJohn1975}. The rest of the methods are mainly based on the maximization of the trace or maximization of the minimum eigenvalue of a matrix.

In the D-optimal framework, maximization of the cost function \eqref{Mm} is given by
\begin{align}
\begin{aligned}%\mbox{\it \small Moment D-optimal problem:~~}
\label{D-opt}
& \underset{\upsilon_o}{\text{minimize}}
& & -\mbox{det}\Big(\mathcal{M}(\upsilon_o)\Big) \\
& \text{subject to}
& & \underline{\upsilon}_o \leq \upsilon_o \leq \overline{\upsilon}_o\\
\end{aligned}
\end{align}
where $\det(\mathcal{M})$ denotes the determinant of $\mathcal{M}$. To be consistent with the optimization literature, the maximization of $z$ is written as the minimization of $-z$ throughout this paper. The solution to this optimization problem will provide the optimal parameters for 2D-MSPC.

For the PC-based image processing, we propose another optimization approach based on the norm of $\mathcal{M}$ as follows:
\begin{align}
\begin{aligned}%\mbox{\it Norm-optimal problem:~~}
\label{norm-opt}
& \underset{\upsilon_o}{\text{minimize}}
& & -\left\| \mathcal{M}(\upsilon_o) \right\| \\
& \text{subject to}
& & \underline{\upsilon}_o \leq \upsilon_o \leq \overline{\upsilon}_o\\
\end{aligned}
\end{align}
where $\| \mathcal{M}\|$ denotes the norm of $\mathcal{M}$. The matrix norm definition is given in Appendix \ref{matrixnorm}.

For $\mu_1=\mu_2=1$, $\mathcal{M}=M+  m= \textstyle \sum_o PC_o^2$, and the D-optimal and norm-optimal problems \eqref{D-opt} and \eqref{norm-opt}  are transferred respectively to 
\begin{align}
\begin{aligned}%\mbox{\it \small Moment D-optimal problem:~~}
\label{D-opt2}
& \underset{\upsilon_o,~o=1,\ldots,O}{\text{minimize}}
& & -\mbox{det}\Big(\mbox{PC}_o^2(\upsilon_o)\Big) \\
& \text{subject to}
& & \underline{\upsilon}_o \leq \upsilon_o \leq \overline{\upsilon}_o\\
& & & \mathcal{M}=M+m\\
\end{aligned}
\end{align}
and,
\begin{align}
\begin{aligned}%\mbox{\it Norm-optimal problem:~~}
\label{norm-opt2}
& \underset{\upsilon_o,~o=1,\ldots,O}{\text{minimize}}
& & -\left\| \mbox{PC}_o^2(\upsilon_o) \right\| \\
& \text{subject to}
& & \underline{\upsilon}_o \leq \upsilon_o \leq \overline{\upsilon}_o\\
& & & \mathcal{M}=M+m\\
\end{aligned}
\end{align}

%{\color{red} It can be written $-\mbox{det}\Big(\mbox{PC}_o(\upsilon)\Big)^2$ or $-\left\| \mbox{PC}_o(\upsilon) \right\|^2$}.

%\begin{Remark}
%The solution to \eqref{D-opt2} and \eqref{norm-opt2} does not require $\mathcal{M}$ at each iteration of optimization. With $\mathcal{M}=M+m$, the parameter vectors are optimized separately for each orientation.% The solution to \eqref{D-opt2} or \eqref{norm-opt2} results in an optimal vector for each orientation.  \hfill{$\square$}
%\end{Remark}

Differences between the optimization sets \{\eqref{D-opt} and \eqref{norm-opt}\} and  \{\eqref{D-opt2} and \eqref{norm-opt2}\} are the followings.
\begin{itemize}
\item[-] In order to solve \eqref{D-opt} or \eqref{norm-opt}, at each iteration of the optimization, the PC matrices of all orientations and $\mathcal{M}$ must be computed. Solving \eqref{D-opt2} and \eqref{norm-opt2} does not require $\mathcal{M}$ at each iteration of the optimization. 
\item[-] The solution to \eqref{D-opt} or \eqref{norm-opt} results in one optimal parameter vector for all orientations. The solution to \eqref{D-opt2} or \eqref{norm-opt2} results in $O$ numbers of optimal parameter vectors (one optimal vector for each orientation).
\end{itemize}

%\begin{Remark}
%The solution to \eqref{D-opt} or \eqref{norm-opt} results in an optimal parameter vector for all orientation. The solution to \eqref{D-opt2} or \eqref{norm-opt2} results in an optimal vector for each orientation.
%\end{Remark}
%>>>>>>>>>>>>>>>>>>>>>>>>>>>>>>>>>>>>>
\subsection{A Sub-Optimal Solution to the Norm-Optimal Problem}
In order to avoid the computation of $\mathcal{M}$ at each iteration of the optimization, a sub-optimal criterion of \eqref{norm-opt} is derived, which optimizes the parameter vector separately at each orientation. Note that the results of this subsection are valid if the matrix norm is consistent and $\mu_i \neq 0$ and $\mu_i \neq 1$ , $i=1,2$.  % {\color{red} The performance of optimal and sub-optimal problems is compared in the terms of accuracy and speed, etc.}  %The following definition and lemmas are required to derive the sub-optimal problem. 

Based on the literature on optimal designs \cite{Fredov1972, Goodwin1977, Atkinson1982, Pronzato1985, Berger1994, Pukelsheim1993}, the covariance matrix satisfies the following inequality
\[
\mbox{cov}(\upsilon) \geq \mathcal{F}^{-1}(\upsilon)
\]
where, $\mbox{cov}$ denotes the covariance matrix, $\upsilon$ is the parameter vector and $\mathcal{F}$ denotes the Fisher information matrix. In this case, $\mathcal{F}^{-1}$ is the lower bound of the covariance of estimations (known as the Cram\' er-Rao bound \cite{Cramer1946, Rao1945}). It is a common approach that in order to reduce the covariance of estimations, the the Cram\' er-Rao bound is reduced by maximizing the Fisher information matrix. 

In this section, we obtain an upper bound of the norm of $M$ and $m$. We then derive a sub-optimal solution, which is based on the maximization of the upper bound of $M$ and $m$. First, we need the matrix norm to be consistent. 
\begin{Definition}(\cite{Lyche2012}, \S8.1.1, p. 178)
A matrix norm is {\em consistent} if it is defined on $\C^{q\times r}$ for all $q,r\in \N$ and the sub-multiplicative property  
\begin{align}
\label{ableqab}
\| AB \| \leq \|A\| ~\|B\|
\end{align}
holds for all matrices $A$ and $B$ for which the product $AB$ is defined. 
%\hfill{$\square$}
\end{Definition}

In this paper, the consistent norm is denoted by a subscript, $\|.\|_c$.

\begin{Lemma} (\cite{Lyche2012}, \cite{Schatten1960})
\label{lemakf}
The Frobenius and all Schatten $p-$norms are sub-multiplicative, i.e.,
\begin{align}
\label{ableqab}
\| AB \|_F &\leq \|A\|_F ~\|B\|_F\\
\| AB \|_p & \leq \|A\|_p ~\|B\|_p
\end{align}
\end{Lemma}

\begin{Lemma} 
\label{lempc}
Let $\mbox{PC}_o$ denotes the PC at orientation $o$ and the 2D-MSPC is computed over $O$ orientations. Then, the following relationships hold between the maximum and minimum moments and $\mbox{PC}_o$.
\begin{align}
M&=\frac{1}{2}\Big(\textstyle \sum_{o} \mbox{PC}_o^2 +\sqrt{\textstyle \sum_{o} \mbox{PC}_o^4 }\Big)\\
m&=\frac{1}{2}\Big(\textstyle \sum_{o} \mbox{PC}_o^2-\sqrt{\textstyle \sum_{o} \mbox{PC}_o^4 }\Big)
\end{align}
{\em Proof}: Appendix \ref{p-lempc}. \hfill{$\square$}
\end{Lemma}

\begin{Proposition}
\label{thpcbnd}
Let $\mbox{PC}_o$ denotes the PC at orientation $o$. The 2D-MSPC is computed over $O$ orientations. If a consistent matrix norm is used, then
\begin{align}
\label{Mbnd}
 &\left\| M \right\|_c  \leq  \textstyle \sum_{o} \left\| \mbox{PC}_o^2 \right\|_c\\
 \label{mbnd}
 &\left\| m \right\|_c  \leq  \textstyle \sum_{o} \left\| \mbox{PC}_o^2 \right\|_c
\end{align}
where, $c$ can take $F$ or $p$. If the subscript $c$ is replaced with $F$, the $F-$norm is applied. If the subscript $c$ is replaced with $p$, the Schatten $p-$norm is applied.
\newline
{\em Proof}: Appendix \ref{p-pcbnd}. \hfill{$\square$}
\end{Proposition}
%{\color{red} It can be written  $\left\| \mbox{PC}_o \right\|^2$}.

\begin{Proposition}
\label{thsubopt}
A sub-optimal solution to the optimization problem \eqref{norm-opt} with the consistent norm is obtained by solving the following optimization problem.
\begin{align}
\begin{aligned}%\mbox{\it Norm-optimal problem:~~}
\label{norm-subopt}
& \underset{\upsilon_o,~o=1,\ldots,O}{\text{minimize}}
& & -(|\mu_1|+|\mu_2|)\left\| \mbox{PC}_o^2(\upsilon_o) \right\|_c \\
& \text{subject to}
& & \underline{\upsilon}_o \leq \upsilon_o \leq \overline{\upsilon}_o\\
\end{aligned}
\end{align}
{\em Proof}: Appendix \ref{p-suboptlem}. \hfill{$\square$}

\end{Proposition}

The differences between \eqref{norm-opt} and \eqref{norm-subopt} are as follows:
\begin{itemize}
\item[-] The solution to \eqref{norm-subopt} is obtained by  the parameter vector optimization individually at each orientation, while the solution to \eqref{norm-opt} is a parameter vector optimal for all orientations. 
\item[-] The optimization problem \eqref{norm-subopt} does not require the computation of $\mathcal{M}$ at each iteration of optimization. To solve \eqref{norm-opt}, the cost function  $\mathcal{M}$ is computed at each iteration of optimization.
\end{itemize} 

The differences between \eqref{norm-opt2} and \eqref{norm-subopt} are as follows:
\begin{itemize}
\item[-] The norm in \eqref{norm-opt2} can be of any type, but in \eqref{norm-subopt}, it has to be consistent and satisfies the sub-multiplicative property. 
\item[-] In \eqref{norm-opt2}, $\mathcal{M}=M+m$, but $\mathcal{M}$ can be any combinations of $M$ and $m$ in \eqref{norm-subopt}.% \hfill{$\square$}
\end{itemize}

%\begin{Remark}
%It is noted that the solution to \eqref{norm-opt2} and \eqref{norm-subopt} The differences between \eqref{norm-opt2} and \eqref{norm-subopt} are: 1) the norm type in \eqref{norm-opt2} can be any, but in \eqref{norm-subopt} has to be consistent and satisfies the sub-multiplicative property. 2) In    \eqref{norm-opt2}, $\mathcal{M}=M+m$, but $\mathcal{M}$ can be any combination of $M$ and $m$ in \eqref{norm-subopt}. \hfill{$\square$}
%\end{Remark}

\subsection{Norm Selection}
Any norm type can be used in \eqref{norm-opt} and \eqref{norm-opt2}. However, if $\mu_i \neq 0$ and $\mu_i \neq 1$, $i=1,2$, only \eqref{norm-subopt} can be applied, not \eqref{norm-opt2}. The prerequisite to apply   \eqref{norm-subopt} is that the sub-multiplicative property must hold for the norm.   

There are also some other features associated with the norm definitions that might be important. If 1-  or $\infty$-norms of the matrix is chosen, because they are associated with the maximum absolute column and row of the matrix, certain features, which are not along that row or column, might be overlooked. The 2-norm, Frobenius-norm ($F-$norm), or Schatten $p-$norm might be more efficient, because they are based on the eigenvalues and singular values, which correspond to the energy of a signal. The 2-norm is the maximum singular value of a matrix: $\|A\|_2=\sigma_{\max}$ \cite{Lyche2012}. The $F-$norm is the square root of the summation of the singular values: $\|A\|_F=(\textstyle \sum_i \sigma_i^2)^{1/2}$ \cite{Lyche2012}, and the Schatten $p-$norm is $\|A\|_p=(\textstyle \sum_i \sigma_i^p)^{1/p}$ \cite{Schatten1960}. Thus, it can be deduced that the 2-norm may amplify the single most-dominant feature, while the $F-$ and $p-$ norms may strengthen features by taking their summation. %Simulation studies show that the norm-based optimization problems are robust to the norm type (see Example 1, \S \ref{sec:ex1}). 

In summary, the 2D-MSPC parameters can be optimally and automatically tuned by maximizing \eqref{Mm}, through \eqref{D-opt} or \eqref{norm-opt}. If $\mathcal{M}=M+m$ is going to be maximized, \eqref{D-opt2} and or \eqref{norm-opt2} can be used. If  $\mu_i \neq 0$ and $\mu_i \neq 1$, $i=1,2$, \eqref{norm-subopt} can provide a sub-optimal solution to \eqref{norm-opt}. This answers Q2.

%\begin{Lemma}(\cite{Meyer2000}, \S5.2, p. 281)
%\label{lem-2norm}
%The 2-norm of a matrix $A$ is the largest singular value of $A$ given by:
%\begin{align}
%\| A\|_2 = \sqrt{\lambda_{\max}(A^{*}A)}.
%\end{align}
%where, $A^{*}$ denotes the conjugate transpose of $A$.\hfill{$\square$}
%\end{Lemma}

\begin{table*}[h]
\caption{Optimal parameters obtained by solving  \eqref{D-opt}, \eqref{norm-opt}, \eqref{D-opt2}, and \eqref{norm-opt2} in Example 1. 
The last column shows the value of the cost function with the optimal parameter vector. }
   \centering
     \begin{tabular}{ l  l  l  l  l}
\hline \hline
		%\multicolumn{2}{c}{Parameters} &  Hints  \cite{Kovesi1996} \\ \hline
       \multirow{1}{*}{Optimization method } & Orientation &   $c_o^*$ &  $g_o^*$ & Cost function opt. value\\ \hline \hline
       \multirow{1}{*}{\eqref{D-opt}}& $o=1,\dots,6$& 0.1212 &    45.4827 & $\mbox{det}(\mathcal{M})^*=7.0088\times 10^{-55}$\\ \hline
                                                                                       
       \multirow{1}{*}{\eqref{norm-opt} based Frobenius norm} &    $o=1,\dots,6$&0.0100      &50.0000  &   $\|\mathcal{M}\|_F^*=114.1408$\\ \hline
						
	\multirow{6}{*}{\eqref{D-opt2}} & $o=1$   & 0.1001  & 49.9997 &  $\mbox{det}(\mbox{PC}_1^2)^*=6.4589\times 10^{-241}$\\ 
	                                                &  $o=2$  & 0.1021   &39.9334 & $\mbox{det}(\mbox{PC}_2^2)^*=1.7131\times 10^{-33}$ \\
						      &  $o=3$  &0.1216    & 46.6111 & $\mbox{det}(\mbox{PC}_3^2)^*=1.4946\times 10^{-43}$\\
						      &  $o=4$  &  0.1946  & 42.9731 & $\mbox{det}(\mbox{PC}_4^2)^*=3.2306\times 10^{-218}$\\
						      &   $o=5$ & 0.1566    &39.2954 & $\mbox{det}(\mbox{PC}_5^2)^*=1.0452\times 10^{-37}$\\
						     &   $o=6$   & 0.1492    &46.3088 & $\mbox{det}(\mbox{PC}_6^2)^*=1.7254\times 10^{-53}$\\% \hline
      \multirow{1}{*}{Average over all orientations} & $o=1,\ldots,6$& 0.1374& 44.1869 & $\mbox{det}(\mathcal{M})^*=\mbox{inf}$ \\  \hline

     \multirow{6}{*}{\eqref{norm-opt2} based Frobenius norm}   &  $o=1$&0.0100    &  50.0000   &  $\|\mbox{PC}_1^2\|_F^*=8.3413\times 10^{3}$\\ 
		                                             & $o=2$ & 0.0100  &   50.0000   &   $\|\mbox{PC}_2^2\|_F^*=6.3400\times 10^{3}$   \\
					           	   &$o=3$ &0.0100    &    50.0000  &   $\|\mbox{PC}_3^2\|_F^*=3.2968\times 10^{3}$   \\
				                  	   &$o=4$	& 0.0100   & 49.9998     &   $\|\mbox{PC}_4^2\|_F^*=2.9181\times 10^{3}$  \\
				                  	   &	$o=5$ &0.0100   &  50.0000   &  $\|\mbox{PC}_5^2\|_F^*=4.7028\times 10^{3}$     \\
				                 	   &	$o=6$ &0.0100 &     50.0000  &  $\|\mbox{PC}_6^2\|_F^*=7.6949\times 10^{3}$\\% \hline
	\multirow{1}{*}{Average over all orientations} & $o=1,\ldots,6$& 0.0100& 50.0000 &$\|\mathcal{M}\|_F^*=3.0379\times 10^{4}$\\ \hline
\end{tabular}
\label{table:resultex1}
\end{table*}

%\section{Visual Convergence}
%\label{sec:visconv}
%The above optimization problems can also be applied to address Q3. Assume that there is one or more parameters in the parameter vector $\upsilon$, in which their optimal values were calculated at their upper or lower limit. For example, assume that $c_o^*=\underline{c}_o$, i.e., the optimal value of $c_o$ is to be $\underline{c}_o$. In such a case, by re-adjusting the limit and solving the corresponding optimization problem, one can assess the visual convergence problem stated Q3. % It is very hard to justify this if the optimal value of a parameter takes a value between the lower and upper limits. 

%>>>>>>>>>>>>>>>>>>>>>>>>>>>>>>>>>>>>>>>>>>>
\section{Illustrative Examples}
\label{sec:ex}

\subsection{Example 1}
\label{sec:ex1}
Consider the Lena's image (Figure \ref{fig:lena-orig}). The objective in this example is to find the optimal values for the parameters of the weighting function while keeping the other parameters constant, where
\begin{align*}
%\label{parvec-ex1}
 & \lambda_{\min o} =3,~ \sigma_{no} =0.55,~ \eta_o=2.1,~ k_o =2,\\& \varepsilon=0.0001 ,~ N_o =4,~ O=6. 
\end{align*}

The parameter vector in this example contains only $c_o$ and $g_o$, thus 
\begin{align*}
\upsilon_o= [ \begin{array}{cc}
    c_o & g_o \end{array}]
    \end{align*}
We set fixed upper and lower limits  $\upsilon_o$ as follows:
\begin{align*}
    & \underline{\upsilon}_o=[ \begin{array}{cc}
    0.1 & 1 \end{array}] \\
    & \overline{\upsilon}_o=[ \begin{array}{cc}
    0.9 & 50 \end{array}].
\end{align*}
With $c_o < 0.1$, almost all spreads of frequencies (even narrow ones) will be kept, but that is not desirable. The change in the slope of the weighting function is also not significant for $g_o > 50$. 

%The proposed methods find the optimal value of the parameters for any range that is determined on $\varepsilon$. 

In order to compare \eqref{D-opt}, \eqref{norm-opt}, \eqref{D-opt2} and \eqref{norm-opt2}, this example focuses on $\mu_1=\mu_2=1$, and therefore,
\begin{align*}
\mathcal{M}=M+m.
\end{align*}

Table \ref{table:resultex1} shows the optimal values of the parameter vector and cost functions. Recall that the solution to \eqref{D-opt} and \eqref{norm-opt} results in one optimal vector for all orientations, and the solution to \eqref{D-opt2} or \eqref{norm-opt2} results in one optimal vector for each orientation. The results of Frobenius-norm optimization are given in Table \ref{table:resultex1}. Recall that the superscript $^*$ represents the optimal value.
% The cost functions of \eqref{D-opt}, \eqref{norm-opt}, \eqref{D-opt2}, and \eqref{norm-opt2} are $-\mbox{det}(\mathcal{M})$, $-\|\mathcal{M}\|_F$, $-\mbox{det}(\mbox{PC}_o^2)$ and $-\|\mbox{PC}_o^2\|_F$, respectively. In the Table the values of $\mbox{det}(\mathcal{M})^*$, $\|\mathcal{M}\|_F^*$, $\mbox{det}(\mbox{PC}_o^2)^*$ and $\|\mbox{PC}_o^2\|_F^*$ are given. 

\begin{figure}[h]
        \begin{subfigure}[b]{0.25\textwidth}
        		\centering 
                \includegraphics[width=\linewidth]{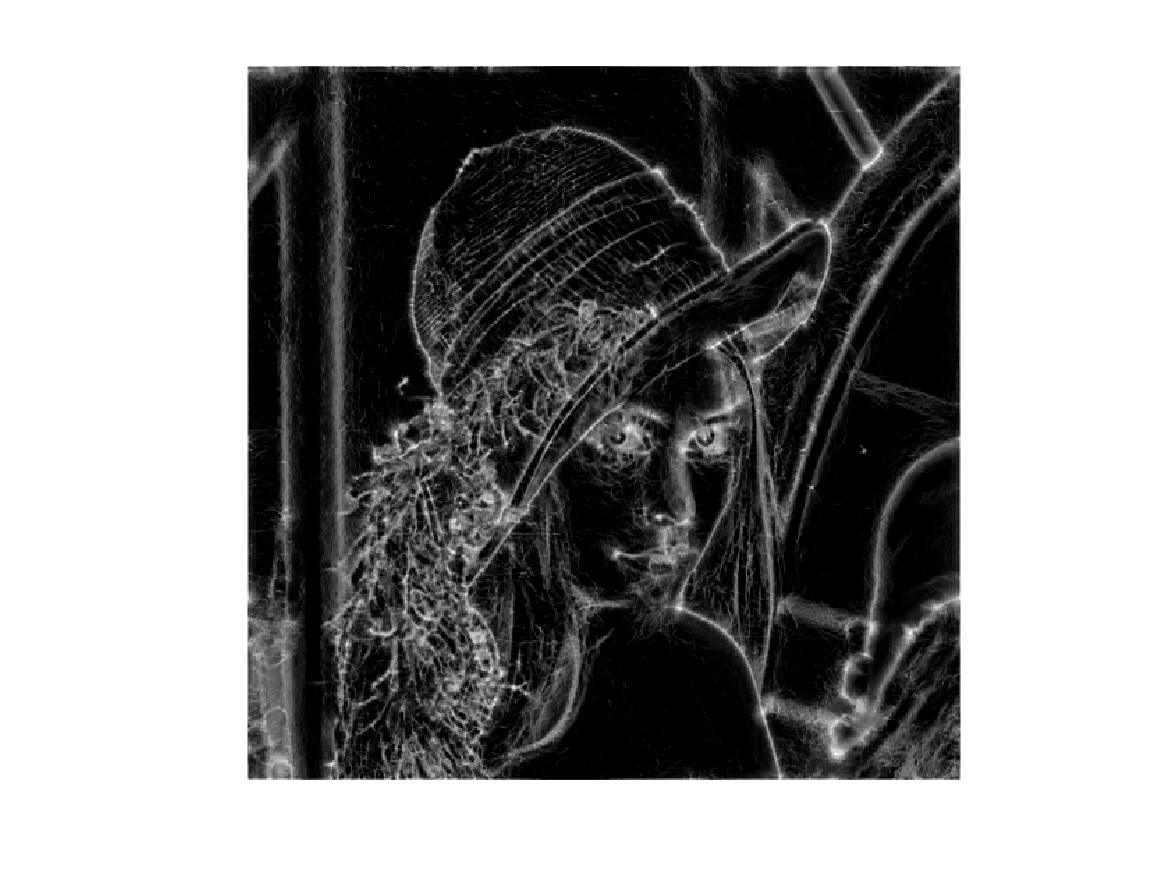}
                \caption{The result of \eqref{D-opt}}
                \label{fig:pc-det-mpm6or}
        \end{subfigure}%
        \begin{subfigure}[b]{0.25\textwidth}
                \includegraphics[width=\linewidth]{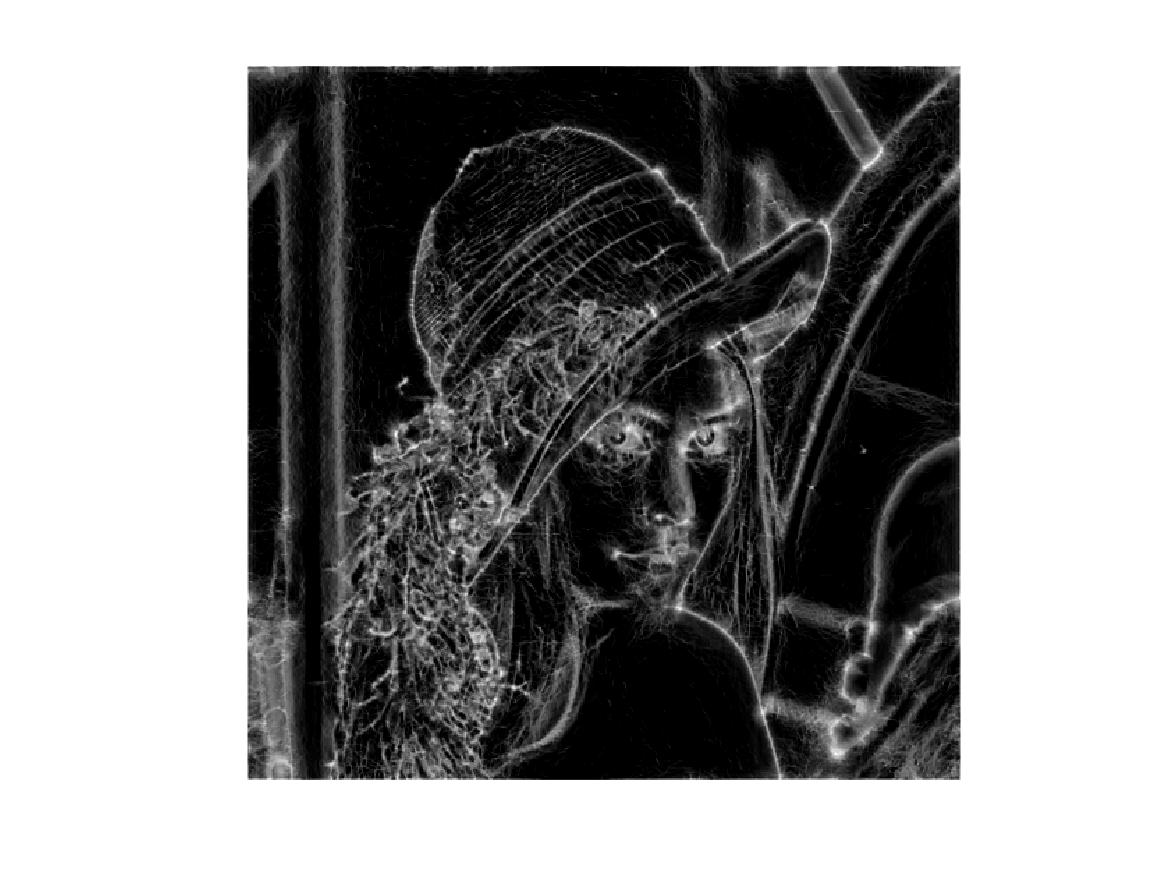}
                \caption{the result of \eqref{norm-opt}}
                \label{fig:pc-fnorm-mpm6or}
        \end{subfigure}\\
        \begin{subfigure}[b]{0.25\textwidth}
                \includegraphics[width=\linewidth]{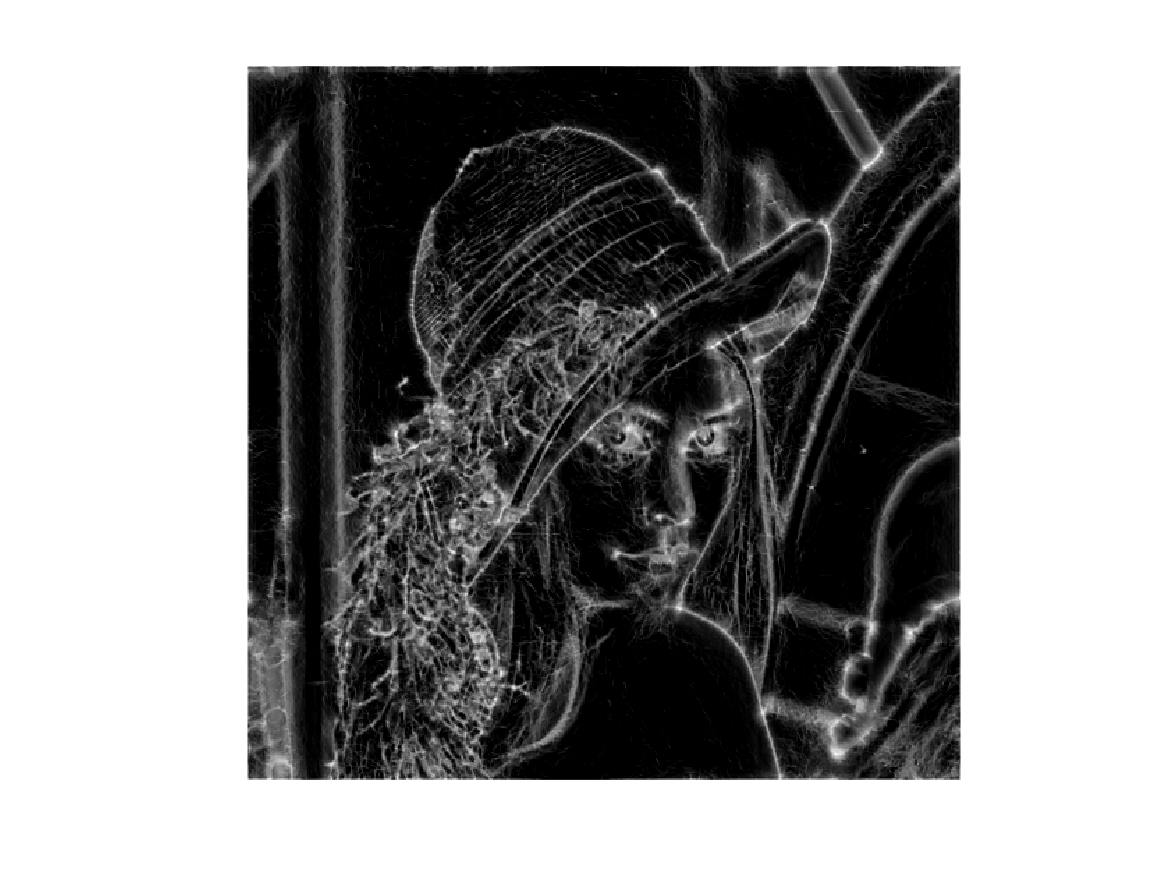}
                \caption{the result of \eqref{D-opt2} }
                \label{fig:pc-det-pc6or}
        \end{subfigure}%
        \begin{subfigure}[b]{0.25\textwidth}
                \includegraphics[width=\linewidth]{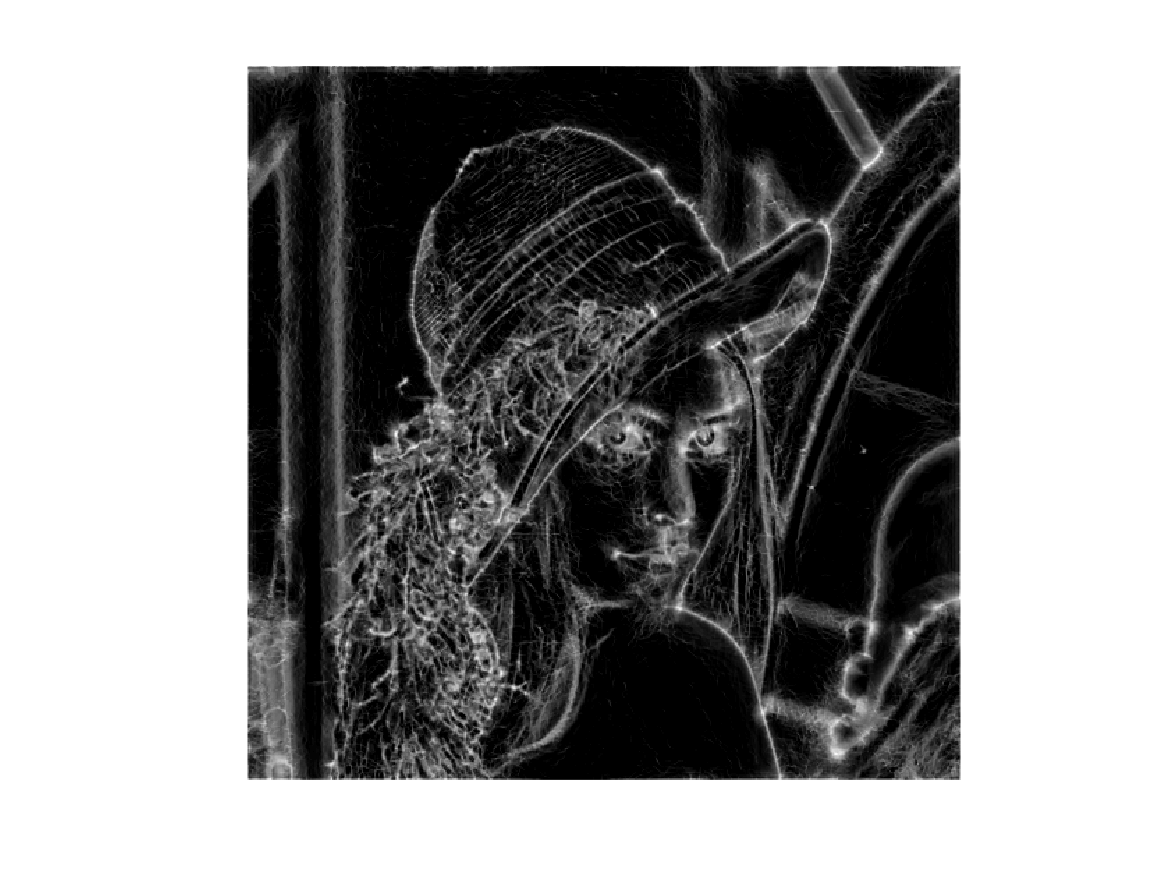}
                \caption{the result of  \eqref{norm-opt2}}
                \label{fig:pc-fnorm-pc6or}
        \end{subfigure}
        \caption{Detection of features from the Lena's image. This is done by optimizing the parameters of the weighting function for 2D-MSPC in Example 1. The image of $\mathcal{M}^*$ is obtained through methods \eqref{D-opt}, \eqref{norm-opt}, \eqref{D-opt2} and \eqref{norm-opt2}.}
         \label{fig:lena-opt}
\end{figure}

Figure \ref{fig:lena-opt} shows the image of $\mathcal{M}^*$ for \eqref{D-opt}, \eqref{norm-opt}, \eqref{D-opt2} and \eqref{norm-opt2}. Features such as lines, edges, Mach bands, corners, and blobs are satisfactorily detected in all images. Visually, it is seen that all four optimization methods result in the same image. 

\begin{figure}[h]
\centering
        \begin{subfigure}[b]{0.25\textwidth}
        		\centering 
                \includegraphics[width=\linewidth]{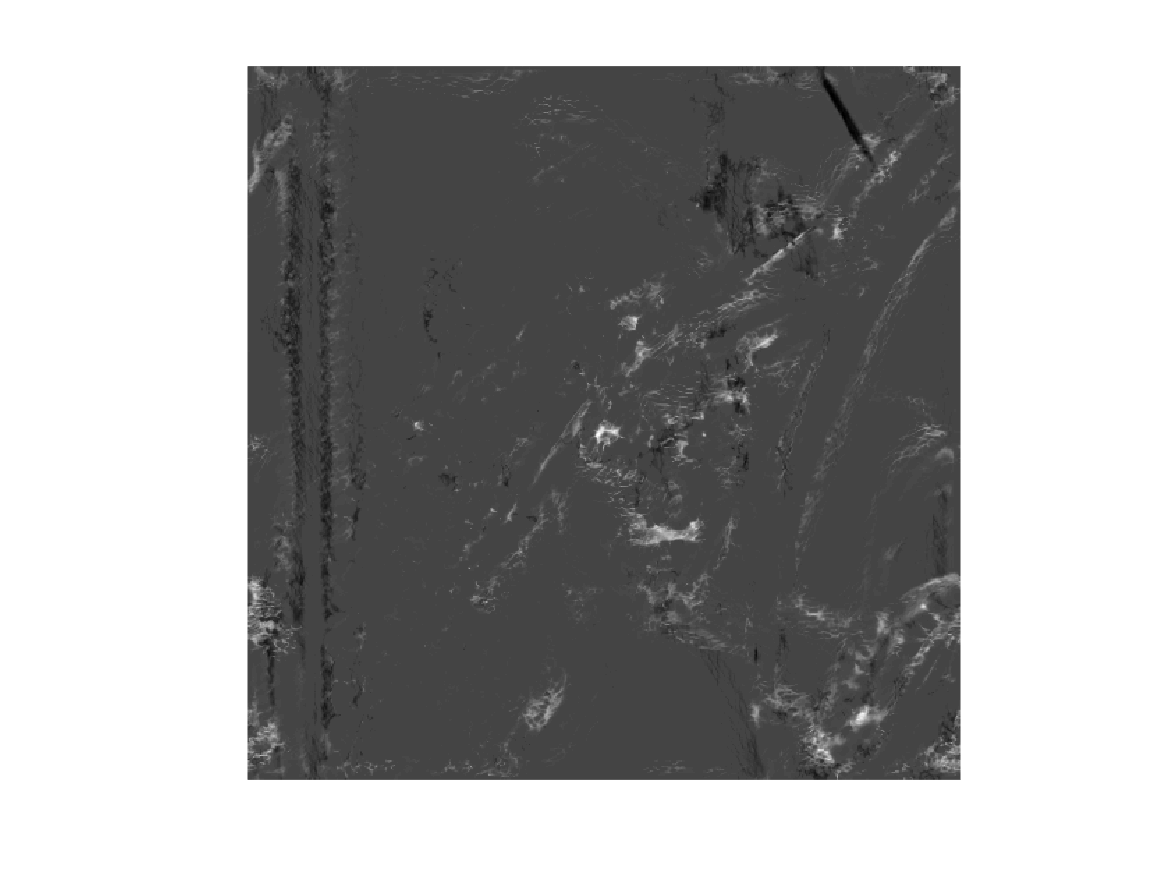}
                \caption{Subtraction Fig. \ref{fig:pc-det-mpm6or} from Fig. \ref{fig:pc-det-pc6or} }
                \label{fig:pc-Fnorm-mpm-full}
        \end{subfigure}\\%
        \begin{subfigure}[b]{0.25\textwidth}
                \includegraphics[width=\linewidth]{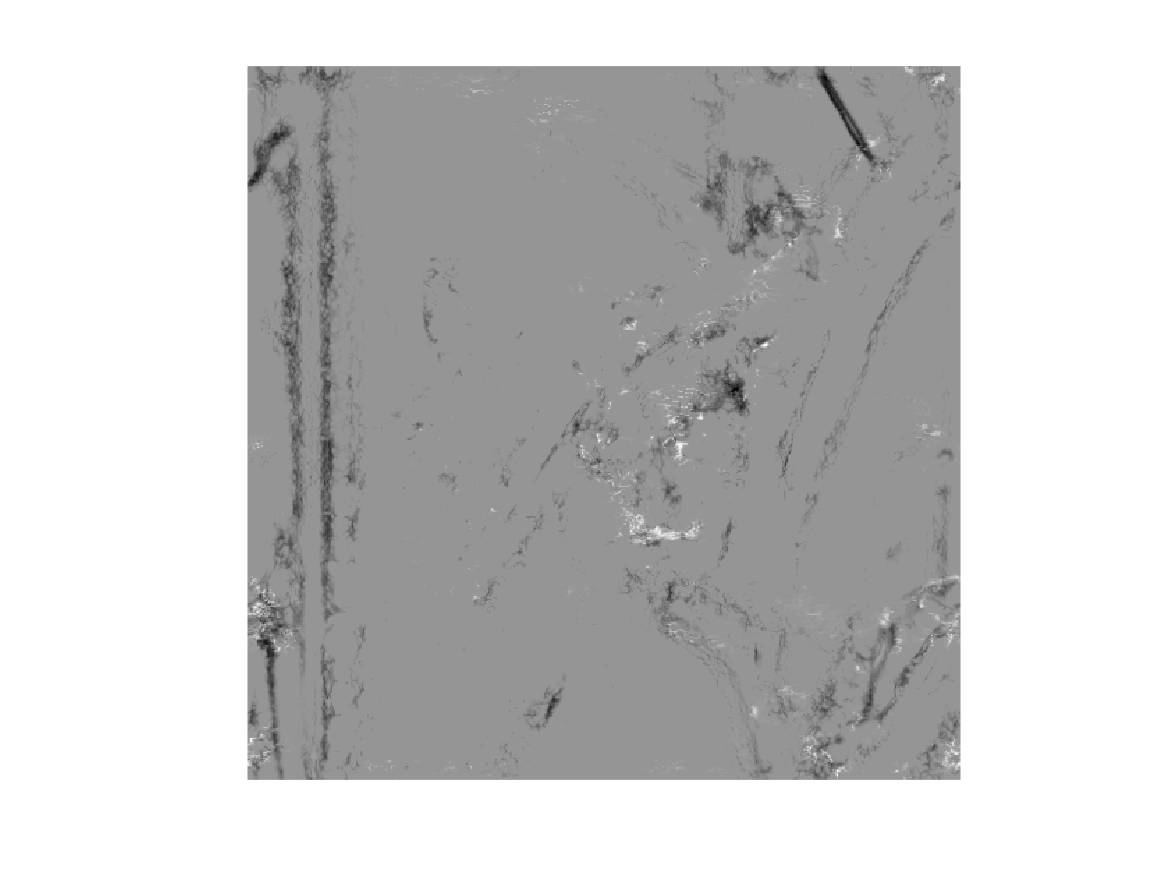}
                \caption{Subtraction Fig. \ref{fig:pc-fnorm-mpm6or} from Fig. \ref{fig:pc-fnorm-pc6or} }
                \label{fig:pc-Fnorm-pc-full}
        \end{subfigure}%\\
        \caption{Comparison between determinant-based optimization methods \eqref{D-opt} and \eqref{D-opt2}, also between norm-based optimization methods \eqref{norm-opt} and \eqref{norm-opt2}.}
         \label{fig:diff-2par}
\end{figure}

Figure \ref{fig:diff-2par} shows the comparison between determinant-based optimization methods \eqref{D-opt} and \eqref{D-opt2}, and norm-based optimization methods \eqref{norm-opt} and \eqref{norm-opt2}. Differences are demonstrated in the subtraction images between the two approaches. Notably, although the determinant and norm optimal frameworks result in a similar performance visually, their subtraction images are not empty, meaning not identical.

%\begin{figure}[h]
%\centering
%\includegraphics[width=0.25\textwidth]{pc_lena_Fnorm_Mpm_m_PC_2par}
%\caption{The subtract between Figures \ref{fig:pc-Fnorm-mpm-full} and \ref{fig:pc-Fnorm-pc-full}}
%\label{fig:diff}
%\end{figure}

\begin{figure}[h]
        \begin{subfigure}[b]{0.25\textwidth}
        		\centering 
                \includegraphics[width=\linewidth]{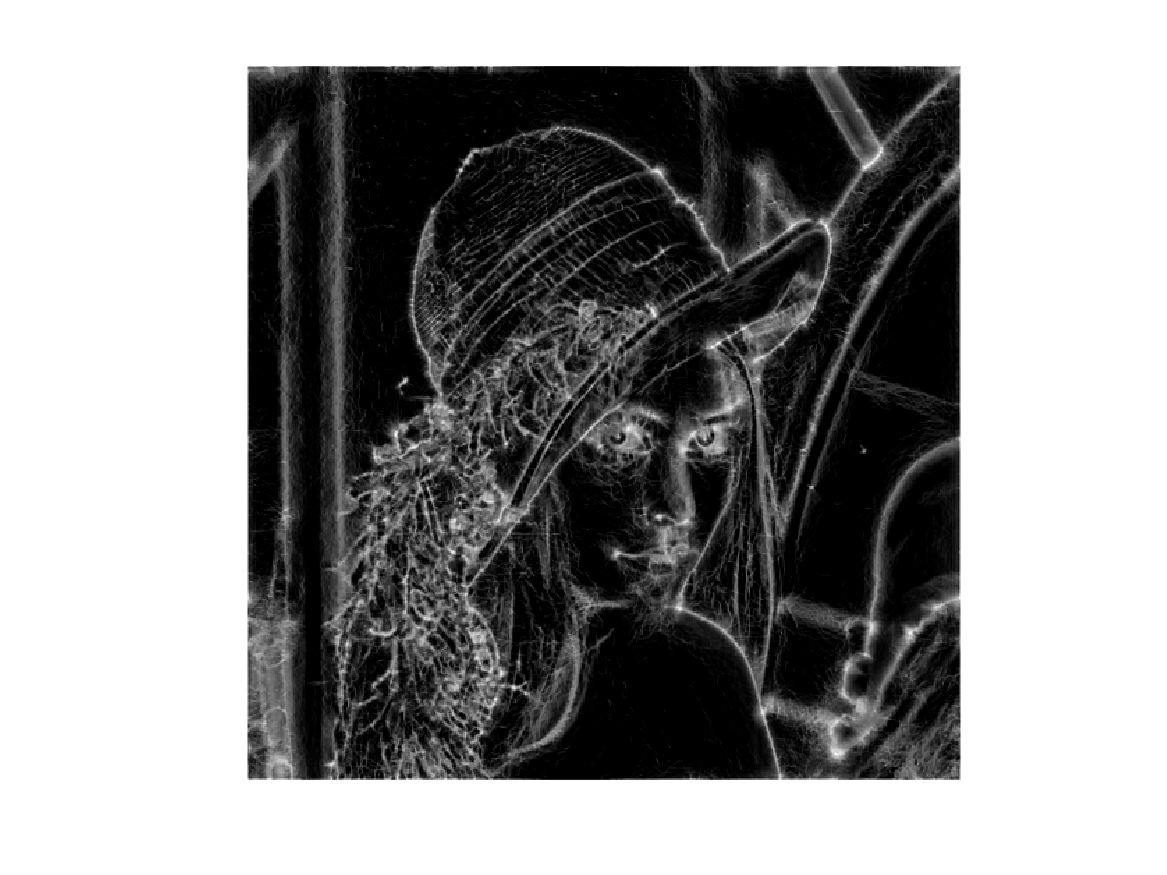}
                \caption{1-norm, $\|\mathcal{M}\|_1^*=219.5310$ }
                \label{fig:pc-1norm-mpm6or}
        \end{subfigure}%
        \begin{subfigure}[b]{0.25\textwidth}
                \includegraphics[width=\linewidth]{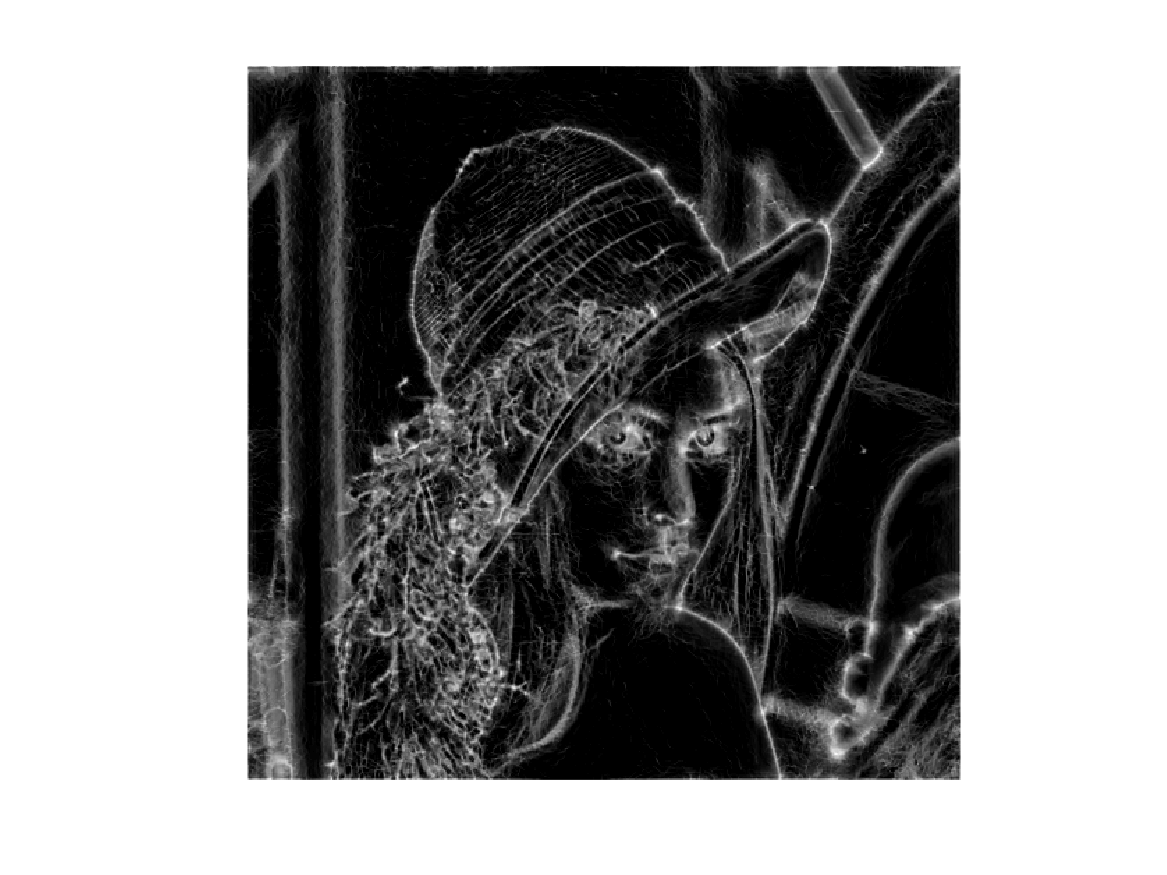}
                \caption{2-norm $\|\mathcal{M}\|_2^*=75.4579$}
                \label{fig:pc-2norm-mpm6or}
        \end{subfigure}\\
        \begin{subfigure}[b]{0.25\textwidth}
                \includegraphics[width=\linewidth]{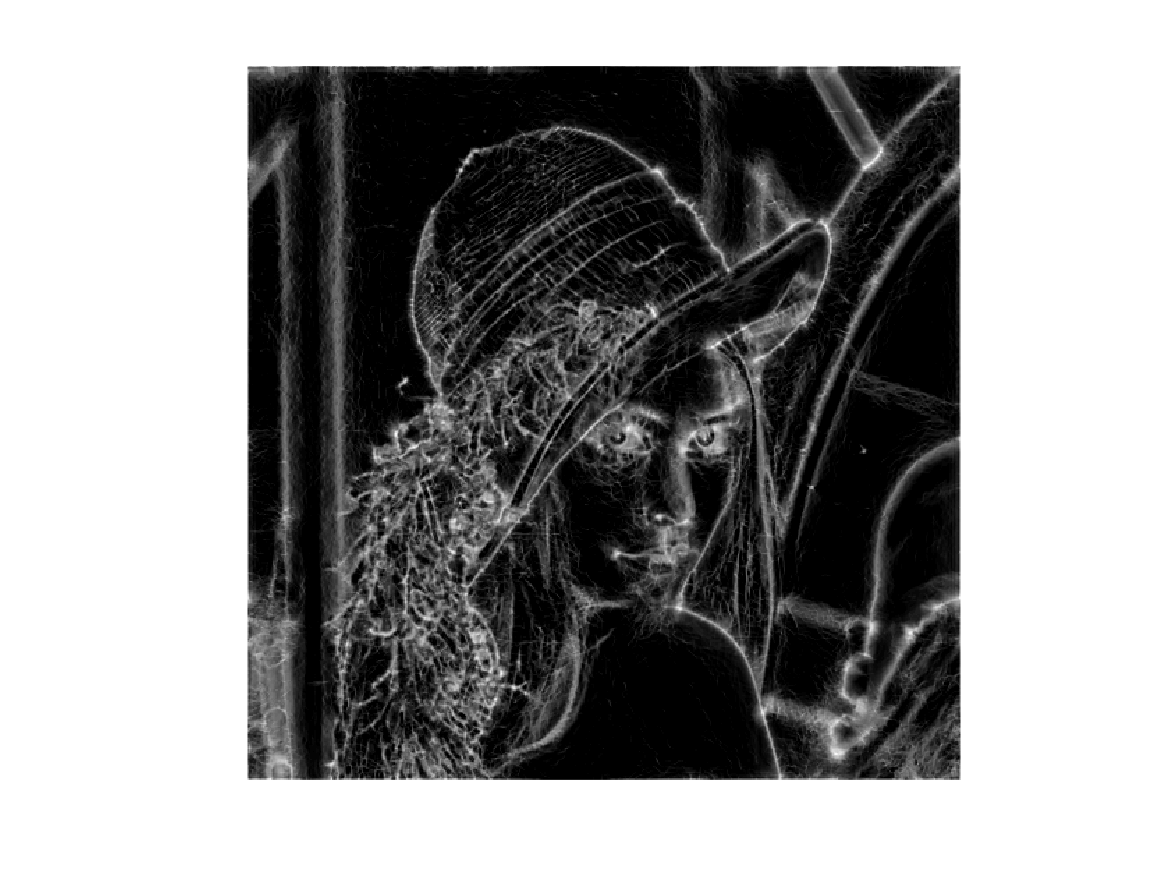}
                \caption{$\infty-$norm $\|\mathcal{M}\|_{\infty}^*=204.6513$ }
                \label{fig:pc-Infnorm-mpm6or}
        \end{subfigure}%
        \begin{subfigure}[b]{0.25\textwidth}
                \includegraphics[width=\linewidth]{pc_lena_Fronorm_Mpm_6or}
                \caption{$F-$norm $\|\mathcal{M}\|_F^*=114.1408$}
                \label{fig:pc-fronormpc6or}
        \end{subfigure}
        \caption{Impact of using different norm types for the detection of features in Lena's image. The image of $\mathcal{M}^*$ and the value of the cost function (calculated for the optimal parameter vector) are achieved from \eqref{norm-opt} based on 1-, 2-, $\infty-$ and $F-$ norms. }
         \label{fig:lena-norm-opt}
\end{figure}

For 1-, 2-, $\infty-$ and $F-$ norms, the optimization method \eqref{norm-opt} results in the optimal parameter vector $\upsilon_o^*=[0.1 ~ 50]$ with $\|\mathcal{M}\|^*_1=219.5310$, $\|\mathcal{M}\|^*_2=75.4579$,  $\|\mathcal{M}\|^*_{\infty}=204.6513$, and $\|\mathcal{M}\|^*_F=114.1408$, respectively. The image of $\mathcal{M}^*$ with optimal settings for the weighting function is shown in Figure \ref{fig:lena-norm-opt}, obtained through \eqref{norm-opt} based on 1-, 2-, $\infty-$ and $F-$ norms. Again, it is seen that all four norm types result in the same visual performance. This example shows that the norm-based optimization of 2D-MSPC for the Lena's image is robust to the norm type.
%The subtract of images is almost empty (not shown here due to page limitations). 

%The optimization problem \eqref{norm-opt2} also results in the optimal parameter vector $\upsilon_o^*=[0.1 ~ 50]$  with different $\|\mbox{PC}_o^2\|^*$ for each orientation and norm type. The images are visually exactly the same as shown in  Figure \ref{fig:lena-norm-opt}. Because of space limitations, both numerical and graphical results are omitted. This example shows that the norm-based optimization of 2D-MSPC for the Lena's image is robust to the norm type.

\begin{Remark}
Because the PC is a matrix with elements between 0 and 1, the determinant is expected to be small. Note that, the larger the image, the smaller the determinant is most likely achieved. Thus, from the numerical point of view, solving such an optimization problem should be more sophisticated than solving the norm-based optimization problem.   \hfill{$\square$}
\end{Remark}

\begin{table*}[h]
\caption{Optimal parameters obtained by solving \eqref{norm-opt} and \eqref{norm-opt2} based on Frobenius norm in Example 2. }
   \centering
     \begin{tabular}{ l  l  l  l  l l l l l}
\hline \hline
		%\multicolumn{2}{c}{Parameters} &  Hints  \cite{Kovesi1996} \\ \hline
       \multirow{1}{*}{Optimization method } & $c_o^*$ &  $g_o^*$ & $\lambda_{\min o}^*$ & $\sigma_{no}^*$ &$\eta_o^*$&$N_o^*$&$O^*$ & Cost function opt. value \\ \hline \hline
       \multirow{1}{*}{\eqref{norm-opt} based $F-$norm}& 0.1 & 49.9994 & 3.055 & 0.4 &3.8315&4& 1& $\|\mathcal{M}\|_F^*=442.0584$ \\ \hline

       \multirow{1}{*}{\eqref{norm-opt2} based $F-$norm}& 0.1 & 50.00 & 2.6894 & 0.4 & 4.0 &4& 1& $\|\mbox{PC}_o^2\|_F^*=6.2033\times 10^{4}$ \\ %\hline
					
	  \multirow{1}{*}{}&  &   &   &  &  && & $\|\mathcal{M}\|_F^*=1.2407\times 10^{5}$ \\ \hline
\end{tabular}
\label{table:resultex2}
\end{table*}

%\begin{figure}[h]
%        \begin{subfigure}[b]{0.25\textwidth}
%        		\centering 
%                \includegraphics[width=\linewidth]{pc_lena_Fnorm_Mpm_full}
%                \caption{The result of \eqref{norm-opt}}
%                \label{fig:pc-Fnorm-mpm-full}
%        \end{subfigure}%
%        \begin{subfigure}[b]{0.25\textwidth}
%                \includegraphics[width=\linewidth]{pc_lena_Fnorm_PC_full}
%                \caption{The result of \eqref{norm-opt2}}
%                \label{fig:pc-Fnorm-pc-full}
%        \end{subfigure}%\\
%        \caption{The results of the optimization of 7 parameters by solving Frobenius norm-based methods \eqref{norm-opt} and \eqref{norm-opt2}. The image of $\mathcal{M}^*$ is shown.  }
%         \label{fig:lena-7par}
%\end{figure}
%
%\begin{figure}[h]
%\centering
%\includegraphics[width=0.25\textwidth]{pc_lena_Fnorm_Mpm_m_PC_full}
%\caption{The subtract between Figures \ref{fig:pc-Fnorm-mpm-full} and \ref{fig:pc-Fnorm-pc-full}}
%\label{fig:diff}
%\end{figure}
%

\begin{figure}[h]
\centering
        \begin{subfigure}[b]{0.25\textwidth}
        		\centering 
                \includegraphics[width=\linewidth]{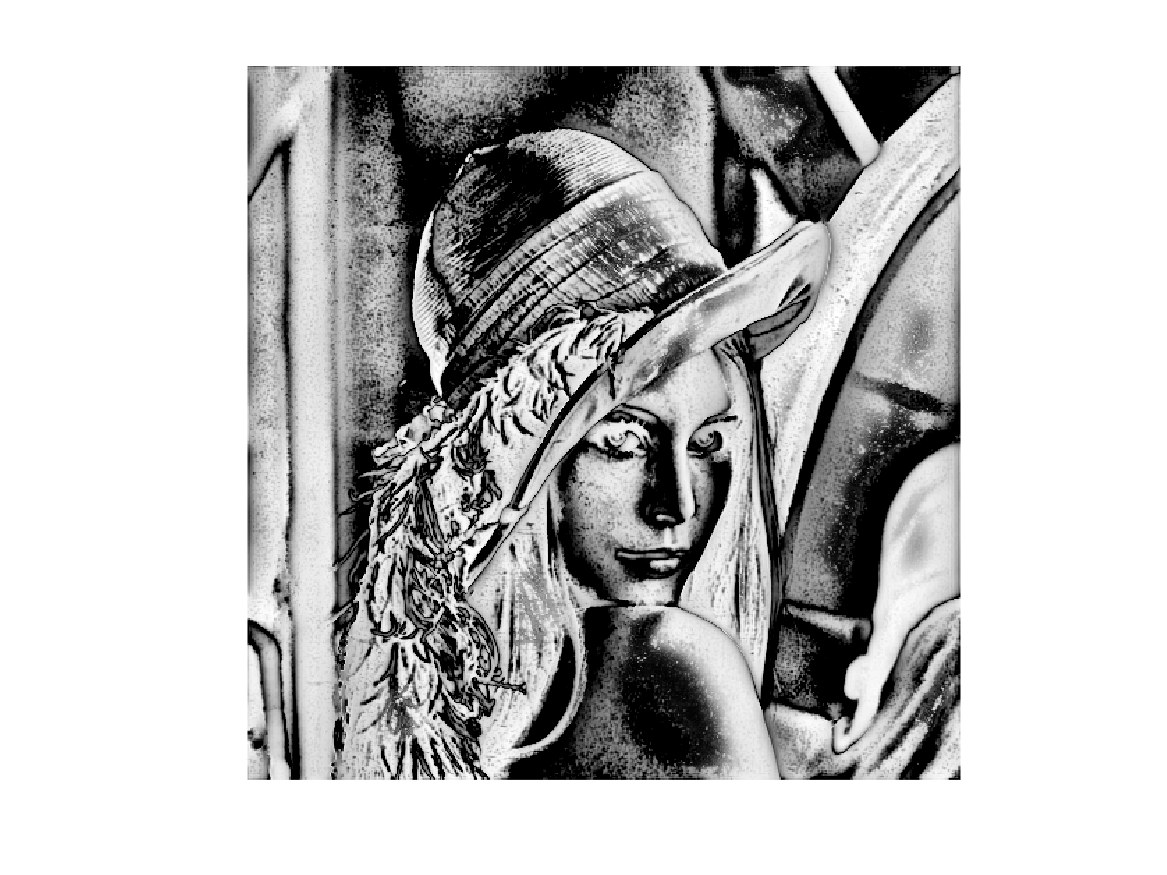}
                \caption{The result of \eqref{norm-opt}}
                \label{fig:pc-Fnorm-mpm-full}
        \end{subfigure}%
        \begin{subfigure}[b]{0.25\textwidth}
                \includegraphics[width=\linewidth]{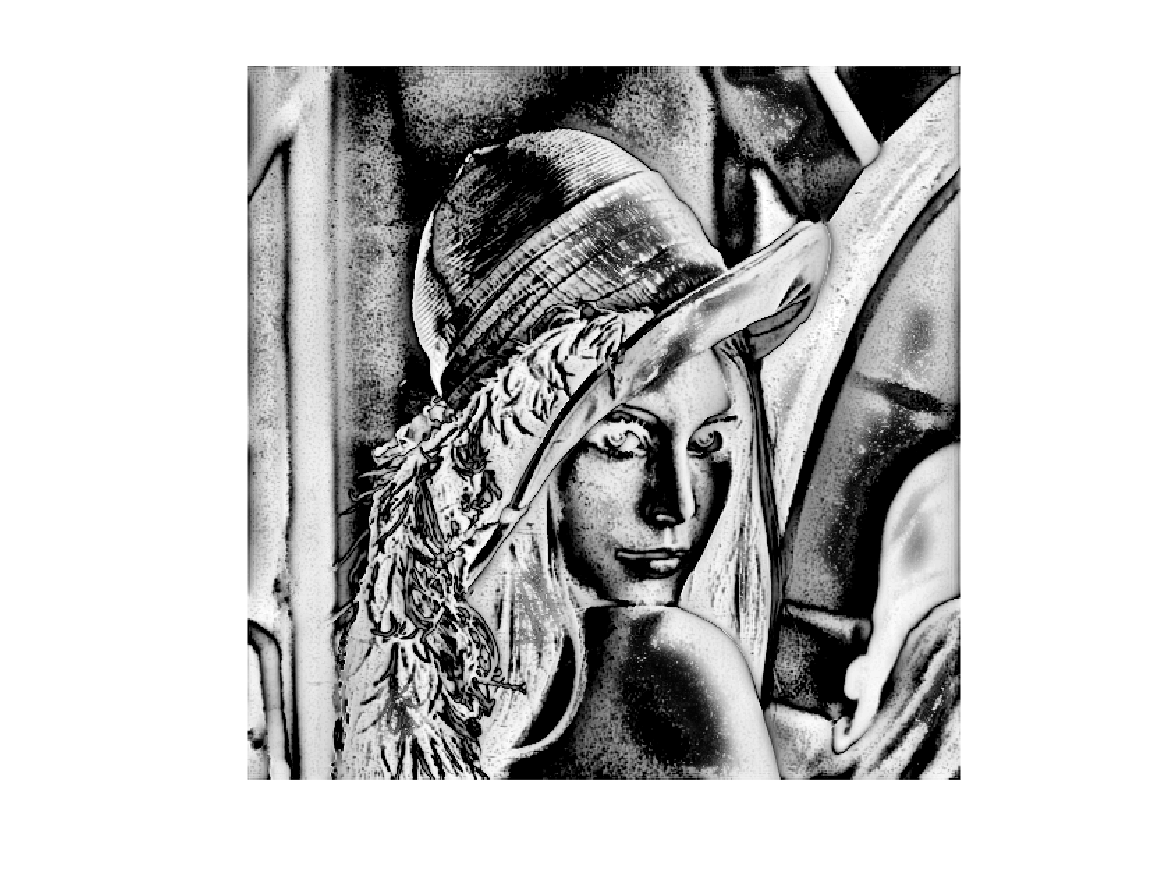}
                \caption{The result of \eqref{norm-opt2}}
                \label{fig:pc-Fnorm-pc-full}
        \end{subfigure}\\
         \begin{subfigure}[b]{0.25\textwidth}
                \includegraphics[width=\linewidth]{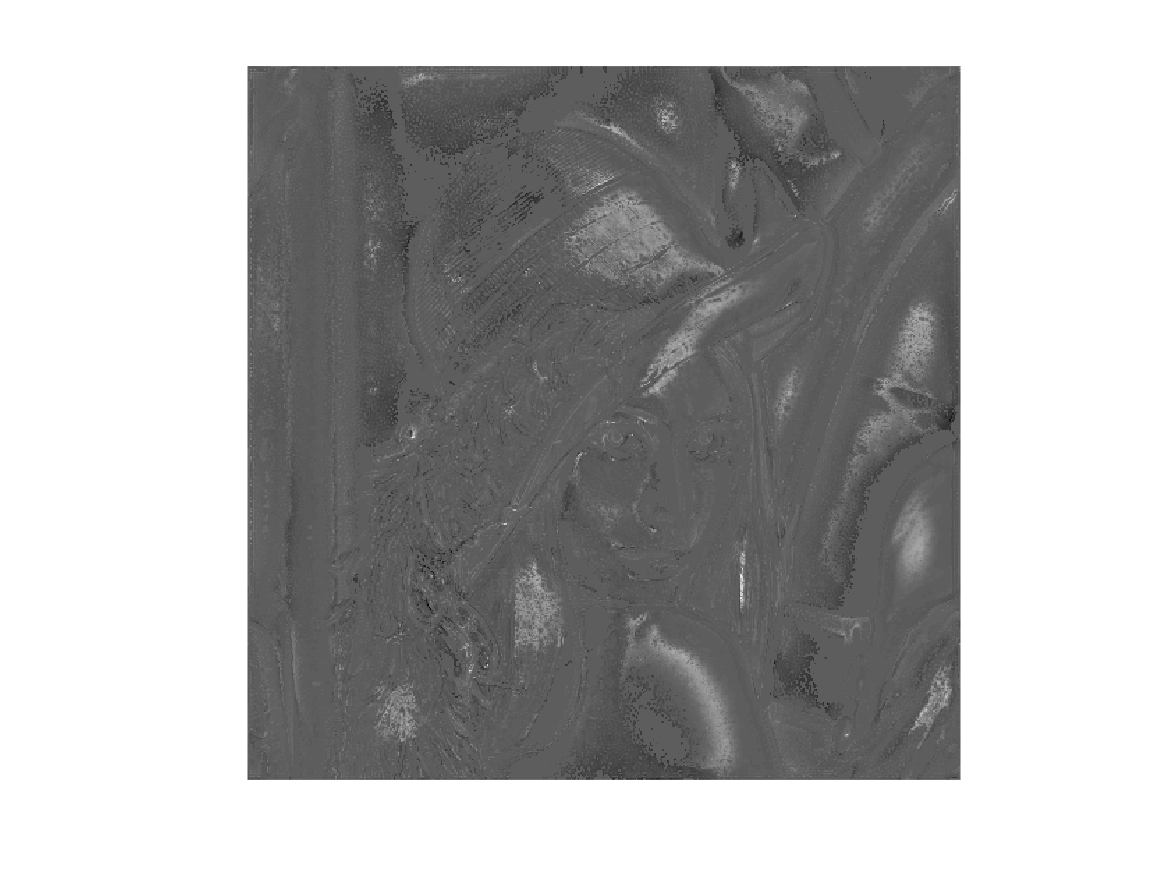}
                \caption{Subtraction Fig. \ref{fig:pc-Fnorm-mpm-full} from Fig. \ref{fig:pc-Fnorm-pc-full}}
                \label{fig:diff}
        \end{subfigure}%\\
        \caption{The results of the optimization of 7 parameters by solving Frobenius norm-based methods in Example 2. Figures \ref{fig:pc-Fnorm-mpm-full} and  \ref{fig:pc-Fnorm-pc-full} shows the  images of $\mathcal{M}^*$ with optimal parameters obtained through \eqref{norm-opt} and \eqref{norm-opt2}, respectively. Figure \ref{fig:diff} shows the difference by subtracting the images.  }
         \label{fig:lena-7par}
\end{figure}

\begin{figure*}[h]
        \begin{subfigure}[b]{0.5\textwidth}
        		\centering 
                \includegraphics[width=1\linewidth]{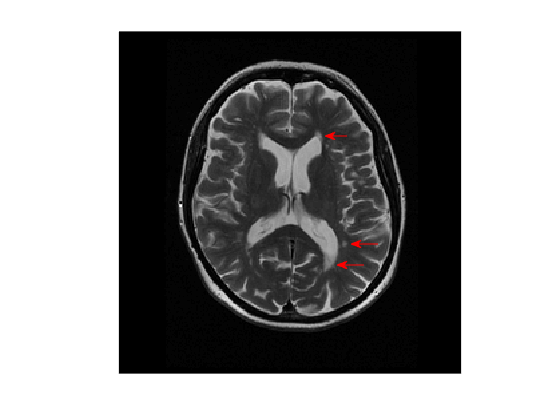}
                \caption{Original MRI}
                \label{fig:IM-0004-0054_patient6_T2-orig}
        \end{subfigure}%
        \begin{subfigure}[b]{0.5\textwidth}
       	       \includegraphics[width=1\linewidth]{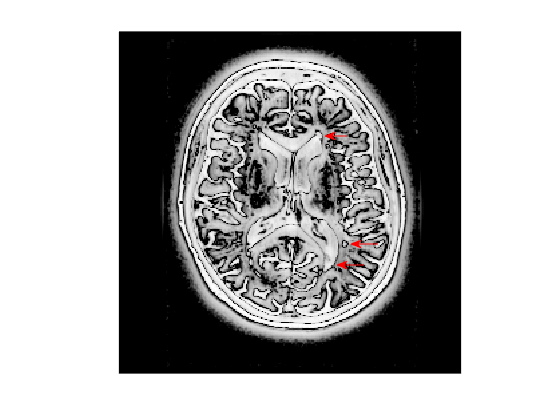}
                \caption{Image of $\mathcal{M}^*$}
                \label{fig:IM-0004-0054_patient6_T2-pc}
        \end{subfigure}%\\
        \caption{Example 3: Detection of MS-MRI lesions using the proposed norm-based 2D-MSPC. Arrows indicate example MS lesions.  The image of $\mathcal{M}^*$ is obtained by using \eqref{norm-opt} with $\mathcal{M}=M$.}
         \label{fig:ms-7par}
\end{figure*}

\subsection{Example 2:}
 The objective of this example is to find the optimal values for the number of scales, number of orientations, and the parameters of the weighting function parameters and bank of filters, with the noise and ill-conditioning parameters fixed to
\begin{align*}
%\label{parvec-ex1}
 k_o =2, \varepsilon=0.0001
\end{align*}

The parameter vector in this example is
\begin{align*}
\upsilon_o= [ \begin{array}{ccccccc}
     c_o &g_o & \lambda_{\min o} & \sigma_{no} & \eta_o & N_o & O \\
     \end{array}]
    \end{align*}
We set fixed upper and lower limits to $\upsilon_o$ as follows
\begin{align*}
&\underline{\upsilon}_o= [ \begin{array}{ccccccc}
    0.1 & 1.0 & 2.0 & 0.4 & 1.0 & 1 & 1 \\
   %  50 & 0.9 & 5 & 1 & 4 & 4 & 6 \\
     \end{array}]\\
&\overline{\upsilon}_o= [ \begin{array}{ccccccc}
  %  1 & 0.1 & 2 & 0.4 & 1 & 1 & 1 \\
     0.9 & 50.0 & 5.0 & 1.0 & 4.0 & 4 & 6 \\
     \end{array}]
    \end{align*}

Because the Lena's IFD has shown robustness to the norm type as demonstrated above, we only test the Frobenius norm-based optimization problems \eqref{norm-opt} and \eqref{norm-opt2}. The resulting optimal values are given in Table \ref{table:resultex2}. Notably, the optimal number of orientations where $O^*=1$ is achieved.  Another noticeable result is that the optimal values for the weighting function parameters do not change, for which the boundary value at the upper and lower limits are found.

%\begin{align*}
%\upsilon_{o}^{*}= [ \begin{array}{ccccccc}
%   49.9994 & 0.1 & 3.055 & 0.4 & 3.8315 & 4 & 1 \\
%     \end{array}]
%    \end{align*}
%with the cost function $\|\mathcal{M}\|^*_F=442.0584$. 
%
%\begin{align*}
%\upsilon_{o}^{*}= [ \begin{array}{ccccccc}
%   50 & 0.1 & 2.6894 & 0.4 & 4 & 4 & 1 \\
%     \end{array}]
%    \end{align*}
%
%$\|\mathcal{M}\|^*_F=1.2407e+05$
%$\|\mbox{PC}_1^2\|^*_F=6.2033e+04$

Figures \ref{fig:pc-Fnorm-mpm-full} and \ref{fig:pc-Fnorm-pc-full} show the images of $\mathcal{M}^*$ with associated optimal parameters $\upsilon_{o}^*$.  By optimizing additional parameters of the 2D-MSPC, it is seen that many features are now better detected than Example 1. This is also consistent with the outcome of cost functions: the value of cost functions is larger than that in Example 1. Figure \ref{fig:diff} shows the difference between the images obtained through \eqref{norm-opt} and \eqref{norm-opt2}.

%By optimizing more parameters of the 2D-MSPC, it is seen that many features are now more detectable compared to the results in Example 1. This is also consistent with the outcome of cost functions: the value of cost functions is larger than that in Example 1. %Examples 1 and 2 confirm that the more the optimal parameters, the better the IFD.

\subsection{Example 3:}
In this example, magnetic resonance imaging (MRI) from the brain of a patient having multiple sclerosis (MS) is studied. MRI plays a key role in diagnosis and management of MS, \cite{Rodriguez2013}. MS is a disease that causes nerve damage in the brain and spinal cord. Characteristically, multi-focal plaques (lesions) can be seen using MRI, showing areas of brightness compared to the surrounding tissue (Figure \ref{fig:IM-0004-0054_patient6_T2-orig}) \cite{Cabezas2014, Schmidt2012, Lorenzo2013}. Our experiments show that maximization of the minimum moment can increase the brightness of brain white matter in MRI. Thus, we chose $\mu_2=0$, and the cost function became $\mathcal{M}=M$. We solved the 2D-MSPC optimization for the $F-$norm with $\underline{\upsilon}_o$ and $\overline{\upsilon}_o$ as given in Example 2. The following optimal parameters were achieved with the optimal cost function $\|\mathcal{M}\|^*_F=234.0708$.  
\begin{align*}
\upsilon_{o}^{*}= [ \begin{array}{ccccccc}
   49.9396 & 0.1001 & 4.7264 & 0.4 & 1.7499 & 3 & 1 \\
     \end{array}]
    \end{align*}
%It is again seen that the optimization led to a single orientation 2D-MSPC computation with $[c_o~ g_o]$ almost at $[\underline{c}_o~ \overline{g}_o]$. 

Figure \ref{fig:IM-0004-0054_patient6_T2-pc} shows the image of $\mathcal{M}$. It is seen that the location and size of lesions are clearly detectable by using the proposed optimal 2D-MSPC.% Image registration and background removal should be applied to the resulting image, which is not the scope of this paper. 

\subsection{Example 4}
%In this example, we show another MRI example from a postmortem brain with MS obtained using a high-field scanner (Figure \ref{fig:NR01196s57R1-orig}). Arrows indicate the MS lesions. MRI pathologies are confirmed by histological analysis. Several regions of interest in the MRI with different severity of tissue damage are highlighted. Figure  \ref{fig:NR01196s57R1-orig} shows the original MRI with regions of interest representing normal appearing white matter (green), diffusely abnormal white matter (yellow) and focal lesions (red). %The objective is too see whether the proposed optimal 2D-MSPC can detect and distinguish these regions. 

\begin{figure*}%[h]
        \begin{subfigure}[b]{0.5\textwidth}
        		\centering 
                \includegraphics[width=1\linewidth]{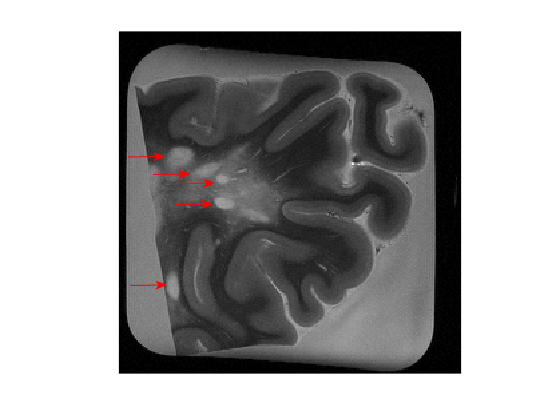}
                \caption{Original Image and regions of interest}
                \label{fig:NR01196s57R1-orig}
        \end{subfigure}%
        \begin{subfigure}[b]{0.5\textwidth}
        \centering 
       	       \includegraphics[width=1\linewidth]{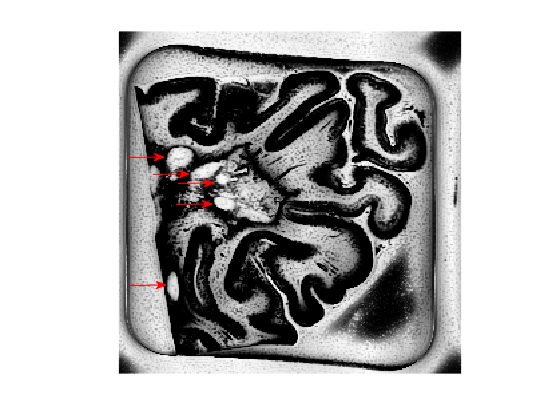}
                \caption{Image of $\mathcal{M}^*$}
                \label{fig:NR01196s57R1-pc}
        \end{subfigure}%\\
%        \begin{subfigure}[b]{0.33\textwidth}
%%        \centering 
%       	       \includegraphics[width=.6\linewidth]{cost_NR01196s57R1_Fnorm_M2pm2_full}
%                \caption{2D-MSPC using \eqref{norm-opt2} \\with $\mathcal{M}=M.*M+m.*m$}
%                \label{fig:ms-pc}
%        \end{subfigure}%\\
        \caption{Example 4: Detection of MS-MRI lesions using the proposed norm-based 2D-MSPC (optimization \eqref{norm-opt}). Arrows indicate the MS lesions. The image of $\mathcal{M}^*$ is obtained by using \eqref{norm-opt} with $\mathcal{M}=M$. }
         \label{fig:NR01196s57R1}
\end{figure*}

In this example, we show another MR image from a postmortem brain with MS obtained using a high-field MR scanner (Figure \ref{fig:NR01196s57R1-orig}). Arrows indicate MS lesions, which have been confirmed by histological analysis.% Figure  \ref{fig:NR01196s57R1-orig} shows the original MRI with regions of interest representing normal appearing white matter (green), diffusely abnormal white matter (yellow) and focal lesions (red). %The objective is too see whether the proposed optimal 2D-MSPC can detect and distinguish these regions. 

Because maximization of the minimum moment increases the brightness of the non-lesion brain areas, we chose $\mu_2=0$. The cost function also became $\mathcal{M}=M$. We solve the 2D-MSPC optimization for the $F-$norm, with $\underline{\upsilon}_o$ and $\overline{\upsilon}_o$ as given in Example 2. The following optimal parameters were achieved with the optimal cost function $\|\mathcal{M}\|^*_F=262.1824$.  
\begin{align*}
\upsilon_{o}^{*}= [ \begin{array}{ccccccc}
   49.9998 & 0.1 & 2 & 0.4 & 3.998 & 4 & 1 \\
     \end{array}]
    \end{align*}
%It is again seen that the optimization led to a single orientation 2D-MSPC computation with with $[c_o~ g_o]$ almost at $[\underline{c}_o ~\overline{g}_o]$. 
Figure \ref{fig:NR01196s57R1-pc} shows the image of $\mathcal{M}^*$. It is seen that the MS lesions are clearly detectable.

%\begin{figure*}%[h]
%        \begin{subfigure}[b]{0.5\textwidth}
%        		\centering 
%                \includegraphics[width=.9\linewidth]{NR03113sA-orig}
%                \caption{Original Image}
%                \label{fig:ms-orig}
%        \end{subfigure}%
%        \begin{subfigure}[b]{0.5\textwidth}
%       	       \includegraphics[width=.9\linewidth]{NR03113sA-pc0}
%                \caption{The result of \eqref{norm-opt2}}
%                \label{fig:ms-pc}
%        \end{subfigure}%\\
%        \caption{The results of optimal 2D-MSPC of a brain magnetic resonance image of a Multiple Sclerosis patient.  The regions plotted by red lines show the focal lesions. The yellow regions denote the dirty appearing white matter. The green region represents the normal appearing white matter.     }
%         \label{fig:ms-7par}
%\end{figure*}

%\begin{figure*}%[h]
%        \begin{subfigure}[b]{0.33\textwidth}
%        		\centering 
%                \includegraphics[width=.6\linewidth]{NR01196s57R1-orig}
%                \caption{Original Image\\ and regions of interest}
%                \label{fig:ms-orig}
%        \end{subfigure}%
%        \begin{subfigure}[b]{0.33\textwidth}
%        \centering 
%       	       \includegraphics[width=.6\linewidth]{cost_NR01196s57R1_Fnorm_M_full}
%                \caption{2D-MSPC using \eqref{norm-opt2}\\ with $\mathcal{M}=M$}
%                \label{fig:ms-pc}
%        \end{subfigure}%\\
%        \begin{subfigure}[b]{0.33\textwidth}
%%        \centering 
%       	       \includegraphics[width=.6\linewidth]{cost_NR01196s57R1_Fnorm_M2pm2_full}
%                \caption{2D-MSPC using \eqref{norm-opt2} \\with $\mathcal{M}=M.*M+m.*m$}
%                \label{fig:ms-pc}
%        \end{subfigure}%\\
%        \caption{The results of optimal 2D-MSPC of a MS MRI. The regions of interest include: normal appearing white matter (shown by green), dirty appearing white matter (shown by yellow) and focal lesions (shown by red). }
%         \label{fig:ms-7par}
%\end{figure*}

%\begin{figure}[h]
%\centering
%\includegraphics[width=0.25\textwidth]{ms-orig}
%\caption{The subtract between Figures \ref{fig:pc-Fnorm-mpm-full} and \ref{fig:pc-Fnorm-pc-full}}
%\label{fig:diff}
%\end{figure}

%>>>>>>>>>>>>>>>>>>>>>>>>>>>>>>>>>>>>>>>>>>>
\section{Conclusions}
In this paper, we have focused on the IFD using a  2D-MSPC method.  2D-MSPC is originally proposed by Peter Kovesi and has shown great potential for the detection of various image features, particularly, lines, edges, corners,  Mach bands, and blobs. However, the parameter setting of 2D-MSPC is typically performed manually based on trial and error, and studies for tuning of such parameters are limited. To enhance the application of this method, we have proposed several optimization frameworks for optimal and automatic tuning of the 2D-MSPC parameters. %We discussed about the visual convergence. 
Through demonstration of several examples including MR images from patients with MS, we show that the ability of IFD can significantly be enhanced.

%This work gives rise to following open problems.
%
%\subsection{Open Problem 1}
%%Throughout the simulations we found out that optimal values of the weighting functions parameters are at the boundaries of upper and lower limits. 
%Let denote the lower and upper limits of $c_o$ and $g_o$ by
%\begin{align*}
%& \underline{c}_o \leq c_o \leq \overline{c}_o\\
%&  \underline{g}_o \leq c_o \leq \overline{g}_o
%\end{align*}
%Is that possible to proof that the optimal $c_o$ and $g_o$ are given by
%\begin{align*}
%c_o^{*}&=\underline{c}_o \mbox{~or~} \overline{c}_o\\
%g_o^{*}&=\underline{g}_o \mbox{~or~} \overline{g}_o
%\end{align*}
%where $c_o^{*}$ and $g_o^{*}$ denote the optimal values of $c_o$ and $g_o$ respectively?
%
%\subsection{Open Problem 2}
%%Throughout the simulations we found out that the optimal value of the parameter vector is the same for all orientations, when norm-based optimization problems are solved. Is that possible to prove it mathematically? 
%Is (or under what conditions is) the optimal value of the parameter vector the same for all orientations, when the norm-based optimization problems are solved?
%
%
%The proof of any of the above problems can significantly expedite the computational process. 
%
%
%

%>>>>>>>>>>>>>>> Appendices
\appendices

\section{Matrix Norm Definition}
\label{matrixnorm}
\begin{Definition}(\cite{Lyche2012}, \S8.1, p.177)
\label{normdef}
A function $\|.\|: \C^{q\times r}$ is called a matrix norm on $ \C^{q\times r}$ if for all $A,B \in  \C^{q\times r}$ and all $c \in \C$
\begin{align}
&\mbox{(Positivity)~}\|A\| \geq 0\\
&\mbox{(homogeneity)~}\| c A\| = |c| ~\|A\| \\
&\mbox{(subadditivity)~} \| A+B\| \leq \|A\|+\|B\|
\end{align}
\end{Definition}

\section{Proof of Lemma \ref{lempc}}
\label{p-lempc}
The replacement of \eqref{alpha}, \eqref{gamma} and \eqref{beta} in \eqref{maxmo} results in
 \begin{align*}
& M=\frac{1}{2}\Big(\textstyle \sum_o \mbox{PC}_o^2 \cos^2(\theta_o)+\sum_o \mbox{PC}_o^2 \sin^2(\theta_o)+\\
     & \sqrt{4\textstyle \sum_o  \mbox{PC}_o^4 \cos^2(\theta_o)\sin^2(\theta_o)+\textstyle \sum_o \mbox{PC}_o^4(\cos^2(\theta_o)-\sin^2(\theta_o))^2}\Big)\\
     &\hspace{5em}=\frac{1}{2}\Big(\textstyle \sum_o \mbox{PC}_o^2 ( \cos^2(\theta_o)+ \sin^2(\theta_o))+\\
   & \hspace{5em}~~ \sqrt{\textstyle \sum_o  \mbox{PC}_o^4 (\cos^2(\theta_o)+\sin^2(\theta_o))^2}\Big)\\&\hspace{5em}=\frac{1}{2}\Big(\textstyle \sum_o \mbox{PC}_o^2+\sqrt{\textstyle \sum_o  \mbox{PC}_o^4}\Big).
\end{align*}

The proof for the minimum moment is exactly the same and omitted.

\section{Proof of Proposition \ref{thpcbnd}}
\label{p-pcbnd}
We need the following Lemmas from linear algebra. 

%\begin{Definition}(\cite{Lyche2012}, \S8.1, p.177)
%\label{normdef}
%A function $\|.\|: \C^{q\times r}$ is called a matrix norm on $ \C^{q\times r}$ if for all $A,B \in  \C^{q\times r}$ and all $c \in \C$
%\begin{align}
%&\mbox{(Positivity)~}\|A\| \geq 0\\
%&\mbox{(homogeneity)~}\| c A\| = |c| ~\|A\| \\
%&\mbox{(subadditivity)~} \| A+B\| \leq \|A\|+\|B\|
%\end{align}
%\end{Definition} 

 \begin{Lemma} 
\label{lemapb}
The following property holds for matrix norms.
\begin{align}
\|A-B\| \leq \|A\| +\|B\|
\end{align}
{\em Proof:} is straightforward using the matrix norm definition \ref{normdef}. 
\end{Lemma}

%\begin{Definition}(\cite{Lyche2012}, \S8.1.1, p. 178)
%A matrix norm is {\em consistent} if it is defined on $\C^{q\times r}$ for all $q,r\in \N$ and the following inequality holds for all matrices $A$ and $B$ for which the product $AB$ is defined.
%\begin{align}
%\label{ableqab}
%\| AB \| \leq \|A\| ~\|B\|
%\end{align}
%holds for all $A,B \in \C^{q\times q}$. %\hfill{$\square$}
%\end{Definition}

\begin{Lemma} (\cite{Lyche2012}, \S8.1.1, p.179)
\label{lemak}
For a consistent matrix norm $\|.\|_c$ on $\C^{q  \times q}$, the following inequality holds 
\begin{align}
\| A^k\|_c \leq \|A\|_c^k, \mbox{~for~} k\in \N.
\end{align}
%\hfill{$\square$}
\end{Lemma}

%\hfill{$\square$}
%\end{Lemma}

The proof of Proposition \ref{thpcbnd} is as follows. By using Lemma \ref{lempc}, we have
\[
\|M\|_c=\frac{1}{2}\left\| \mbox{PC}_1^2+\ldots+\mbox{PC}_O^2+\sqrt{\mbox{PC}_1^4+\ldots+\mbox{PC}_O^4}\right\|_c
\]
From Lemma \ref{lemapb}, thus
\[
\|M\|_c \leq \frac{1}{2}\Big(\left\| \mbox{PC}_1^2\right\|_c+\ldots+\left\|\mbox{PC}_O^2\right\|_c+\left\|\sqrt{\mbox{PC}_1^4+\ldots+\mbox{PC}_O^4}\right\|_c\Big)
\]
By using Lemma \ref{lemak}, it yields
\begin{align*}
\|M\|_c& \leq \frac{1}{2}\Big( \left\| \mbox{PC}_1^2\right\|_c+\ldots+\left\|\mbox{PC}_O^2\right\|_c+\sqrt{\left\|\mbox{PC}_1^4+\ldots+\mbox{PC}_O^4\right\|_c}\Big)\\
& \leq \frac{1}{2}\Big(\left\| \mbox{PC}_1^2\right\|_c+\ldots+\left\|\mbox{PC}_O^2\right\|_c+\sqrt{\left\|\mbox{PC}_1^4\right\|_c+\ldots+\left\|\mbox{PC}_O^4\right\|_c}\Big)
\end{align*}
Lemma \ref{lemak} is again used, which yields
\begin{align*}
&\|M\|_c \leq \\
&\frac{1}{2}\Big(\left\| \mbox{PC}_1^2\right\|_c+\ldots+\left\|\mbox{PC}_O^2\right\|_c+\sqrt{\left\|\mbox{PC}_1^2\right\|_c^2+\ldots+\left\|\mbox{PC}_O^2\right\|_c^2}\big)
\end{align*}
From algebra, $\sqrt{a_1^2+\ldots+a_O^2}\leq a_1+\ldots+a_O$ for $a_i\geq0$, thus, 
\begin{align*}
\|M\|_c &\leq \frac{1}{2}\Big(\left\| \mbox{PC}_1^2\right\|_c+\ldots+\left\|\mbox{PC}_O^2\right\|_c+\left\|\mbox{PC}_1^2\right\|_c+\ldots+\left\|\mbox{PC}_O^2\right\|_c\Big)\\
&=\frac{1}{2}\Big(2\left\| \mbox{PC}_1^2\right\|_c+\ldots+2\left\|\mbox{PC}_O^2\right\|_c\Big)
\end{align*}
and the proof of \eqref{Mbnd} is complete. 

By using Lemma \ref{lempc}, we have
\[
\|m\|_c=\frac{1}{2}\Big(\left\| \mbox{PC}_1^2+\ldots+\mbox{PC}_O^2-\sqrt{\mbox{PC}_1^4+\ldots+\mbox{PC}_O^4}\right\|_c\Big)
\]
From the matrix norm property, thus
\[
\|m\|_c \leq \frac{1}{2}\Big(\left\| \mbox{PC}_1^2\right\|_c+\ldots+\left\|\mbox{PC}_O^2\right\|_c+\left\|\sqrt{\mbox{PC}_1^4+\ldots+\mbox{PC}_O^4}\right\|_c\Big)
\]
and the rest of the proof is exactly the same as above. 

\section{Proof of Proposition \ref{thsubopt}}
\label{p-suboptlem}
From the norm definition, \eqref{Mm}, \eqref{Mbnd} and \eqref{mbnd},  it is deduced that 
\begin{align*}
 % \label{Mmbnd}
\left\| \mathcal{M}\right\|_c  \leq  (|\mu_1|+|\mu_2|) \textstyle \sum_{o} \left\| \mbox{PC}_o^2 \right\|_c
\end{align*}
A suboptimal solution to \eqref{norm-opt}, with the consistent norm, is obtained by increasing the upper found of $\left\| \mathcal{M}\right\|_c$. This completes the proof.

\section*{Acknowledgment}
First author would like to thank Shrushrita Sharma and Glen Pridham for some discussions. Funding from the MS Society of Canada, Natural Sciences and Engineering Council of Canada, and Alberta Innovates -- Health Solutions is also acknewledged.

\ifCLASSOPTIONcaptionsoff
  \newpage
\fi

\end{document}